\documentclass[11pt]{article}
\usepackage[margin=1in]{geometry}
\usepackage{amsmath,amssymb,amsthm,mathtools}
\usepackage{booktabs,tabularx,array,microtype}
\usepackage{algorithm,algpseudocode}
\usepackage[numbers,sort&compress]{natbib}
\usepackage{graphicx}
\usepackage{hyperref}
\hypersetup{
 hidelinks,
 pdftitle={Near-Optimal Complexity of Finding Approximate Local Minima of Nonconvex Finite Sums},
 pdfauthor={Qihao Zhou}
}

\newtheorem{theorem}{Theorem}
\newtheorem{lemma}{Lemma}
\newtheorem{proposition}{Proposition}
\newtheorem{corollary}{Corollary}
\theoremstyle{definition}
\newtheorem{assumption}{Assumption}
\newtheorem{definition}{Definition}
\theoremstyle{remark}

\newcommand{\R}{\mathbb R}
\newcommand{\E}{\mathbb E}
\newcommand{\Prob}{\mathbb P}
\newcommand{\norm}[1]{\lVert#1\rVert}
\newcommand{\abs}[1]{\lvert#1\rvert}
\newcommand{\eps}{\varepsilon}
\newcommand{\tO}{\widetilde{O}}
\newcommand{\tTheta}{\widetilde{\Theta}}
\newcommand{\supp}{\operatorname{supp}}
\newcommand{\Lip}{\operatorname{Lip}}
\newcommand{\spn}{\operatorname{span}}
\newcommand{\sign}{\operatorname{sign}}
\newcommand{\Fas}{\mathcal F_{\mathrm{AS}}}
\newcommand{\Fis}{\mathcal F_{\mathrm{IS}}}
\newcommand{\N}{\mathfrak N}
\newcommand{\Nzr}{\mathfrak N^{\mathrm{zr}}}
\newcolumntype{Y}{>{\raggedright\arraybackslash}X}
\usepackage{titlesec}
\usepackage{needspace}
\titlespacing*{\section}{0pt}{2.8ex plus .6ex minus .2ex}{1.3ex plus .2ex}
\titlespacing*{\subsection}{0pt}{2.2ex plus .5ex minus .2ex}{1.0ex plus .2ex}
\titlespacing*{\paragraph}{0pt}{1.6ex plus .3ex minus .2ex}{1em}
\AtBeginDocument{%
 \setlength{\abovedisplayskip}{10pt plus 2pt minus 2pt}%
 \setlength{\belowdisplayskip}{10pt plus 2pt minus 2pt}%
 \setlength{\abovedisplayshortskip}{6pt plus 2pt minus 1pt}%
 \setlength{\belowdisplayshortskip}{6pt plus 2pt minus 1pt}%
 \setlength{\jot}{4pt}%
}

\title{Near-Optimal Complexity of Finding Approximate\\Local Minima of Nonconvex Finite Sums}
\author{
  Qihao Zhou\thanks{Independent Researcher. Email: \texttt{qihao\_zhou@berkeley.edu}.}
}
\date{}

\begin{document}
\maketitle
\begin{abstract}
We study the component-oracle complexity of finding approximate local minima of nonconvex finite sums under mean-squared averaged smoothness with parameter $L$ and a $\rho$-Lipschitz average Hessian. Each query returns the value, gradient, and one Hessian--vector product of a component. We show that the complexity in the high-accuracy regime $\eps\le c_0L^2/(n\rho)$ is $\tTheta(n^{3/4}L^{1/2}\rho^{1/4}\Delta\eps^{-7/4})$ for a sufficiently large initial gap $\Delta$. We attain the upper bound by adapting the inverse solver and eigenvector routine of FastCubic to averaged smoothness. Our lower bound assigns the bridge edges of a nested chain to random components; a replicated spectral construction gives the same bound for unrestricted randomized algorithms. Under individual smoothness, a factor of $n^{1/4}$ remains for $\eps\le cL^2/(n^2\rho)$, up to logarithmic factors. In the complementary low-accuracy regime, we propose PathVR, which uses $\tO(n+\sqrt n\,L\Delta\eps^{-2})$ queries and matches a lower bound from dense weak hiding.
\end{abstract}

\section{Introduction}\label{sec:intro}
We study the finite-sum problem
\begin{equation}\label{eq:problem}
 \min_{x\in\R^{d}}\ f(x):=\dfrac1n\sum_{i=1}^{n}f_i(x),
\end{equation}
where each component $f_i(\cdot)$ is possibly nonconvex. We seek an approximate local minimum in the sense of \citet{nesterov}: a point $\hat x$ satisfying
\begin{equation}\label{eq:target}
 \norm{\nabla f(\hat x)}\le\eps\qquad\text{and}\qquad
 \lambda_{\min}\bigl(\nabla^2 f(\hat x)\bigr)\ge-\sqrt{\rho\eps},
\end{equation}
where $\rho$ is the Lipschitz constant of the average Hessian $\nabla^2f(\cdot)$. The curvature condition excludes saddle points with sufficiently negative curvature. Each oracle call returns the value, gradient, and one Hessian--vector product of a selected component.

\citet{fastcubic} proposed FastCubic, which approximately solves cubic-regularized Newton subproblems~\cite{nesterov} using accelerated variance reduction and shift-and-invert eigenvector routines. When every component Hessian is bounded by $L$, it finds a point satisfying \eqref{eq:target} with
\begin{equation}\label{eq:fastcubic}
 \tO\!\left(n\rho^{1/2}\Delta\eps^{-3/2}+n^{3/4}L^{1/2}\rho^{1/4}\Delta\eps^{-7/4}\right)
\end{equation}
oracle calls, where $\Delta$ bounds the initial objective gap. The first term accounts for full gradients at the cubic steps. The second comes from accelerated variance-reduced solves of shifted Hessian systems and dominates in the \emph{high-accuracy regime} $\eps\le L^2/(n\rho)$. For $n=1$, it matches the first-order lower bounds for deterministic~\cite{zhousharp,cjz} and randomized~\cite{zhangrand} algorithms. Our averaged-smoothness adaptation uses Katyusha~X~\cite{katx} for these solves. These results leave open whether the finite-sum factor $n^{3/4}$ is necessary.

In this paper, we show that the high-accuracy complexity is $\tTheta(n^{3/4}L^{1/2}\rho^{1/4}\Delta\eps^{-7/4})$ under \emph{mean-squared averaged smoothness}, a weaker condition than the individual smoothness used by FastCubic. The upper bound follows by adapting its inverse solver and eigenvector routine to a second-moment bound on the component Hessians.

Our lower bound embeds the nested chain of \citet{cjz} into a finite sum. Following \citet{zhouncsc}, we assign each quadratic bridge edge to a random component. Discovering an edge then costs order $n$ queries, while averaged smoothness permits bridges long enough to give the factor $n^{3/4}$ (Theorem~\ref{thm:nested}). Replicating the spectral construction of \citet{zhangrand} extends the bound to unrestricted randomized algorithms (Theorem~\ref{thm:highlower}). Under individual smoothness, the bounds differ by at most a factor of $n^{1/4}$ up to logarithmic factors for $\eps\le cL^2/(n^2\rho)$.

In the complementary low-accuracy regime, we propose PathVR (path-length-based variance reduction), which uses $\tO(n+\sqrt n\,L\Delta\eps^{-2})$ oracle calls. The dense-weak-hiding instances of \citet{densehiding} give a matching lower bound (Theorem~\ref{thm:low}).

\section{Preliminaries and Related Work}\label{sec:prelim}
We use $\norm{\cdot}$ for the Euclidean norm of vectors and the spectral norm of matrices. We write $[q]=\{1,\ldots,q\}$, $e_k$ for the $k$-th standard basis vector, and $\supp(v)$ for the set of nonzero coordinates of $v$. The notation $O(\cdot)$, $\Omega(\cdot)$, and $\Theta(\cdot)$ hides universal constants; $\tO(\cdot)$ and $\tTheta(\cdot)$ also hide logarithmic factors in the problem parameters, accuracy, failure probability, and dimension.

\subsection{Function Classes and Approximate Local Minima}
The components $f_1,\ldots,f_n:\R^d\to\R$ are twice continuously differentiable. We distinguish the following two smoothness conditions.
\begin{assumption}\label{ass:as}
For all $x,x'\in\R^d$,
\begin{equation}\label{eq:as}
 \dfrac1n\sum_{i=1}^{n}\norm{\nabla f_i(x)-\nabla f_i(x')}^2
 \le L^2\norm{x-x'}^2.
\end{equation}
\end{assumption}
\begin{assumption}\label{ass:is}
For all $i\in[n]$ and $x\in\R^d$, we have $\norm{\nabla^2 f_i(x)}\le L$.
\end{assumption}
We refer to Assumptions~\ref{ass:as} and~\ref{ass:is} as averaged smoothness (AS) and individual smoothness (IS), respectively. Individual smoothness implies averaged smoothness, whereas \eqref{eq:as} allows a component Hessian of norm $\sqrt nL$. Letting $x'\to x$ in \eqref{eq:as} gives, with $H_i=\nabla^2f_i(x)$ and $H=\nabla^2f(x)$,
\begin{equation}\label{eq:asmatrix}
 \dfrac1n\sum_{i=1}^nH_i^2\preceq L^2I,\qquad
 \dfrac1n\sum_{i=1}^n(H_i-H)^2\preceq L^2I,\qquad \norm{H}\le L.
\end{equation}
We impose the following condition on the average Hessian and the initial gap.
\begin{assumption}\label{ass:hess}
For all $x,x'\in\R^d$, we have $\norm{\nabla^2f(x)-\nabla^2f(x')}\le\rho\norm{x-x'}$, and the initial point satisfies $f(0)-\inf_xf(x)\le\Delta$.
\end{assumption}
As in \citet{fastcubic}, this condition concerns the average Hessian. The component Hessian-Lipschitz constants may be larger, although our high-accuracy hard instances have individually $\rho$-Hessian-Lipschitz components. We write $\Fas(n,L,\rho,\Delta)$ and $\Fis(n,L,\rho,\Delta)$ for the classes of finite sums in all finite dimensions that satisfy Assumption~\ref{ass:hess} together with Assumption~\ref{ass:as} and Assumption~\ref{ass:is}, respectively. Throughout, we assume $\eps\le L^2/\rho$, since otherwise the bound $\norm{\nabla^2f(x)}\le L$ already implies the curvature condition in \eqref{eq:target}.

\subsection{Oracle and Algorithm Classes}
We count queries to the following component oracle.
\begin{definition}\label{def:oracle}
Given an index $i\in[n]$, a point $x$, and a direction $v$, the component Hessian--vector oracle returns $(f_i(x),\nabla f_i(x),\nabla^2f_i(x)v)$. Each call counts as one query, so a full gradient or a full Hessian--vector product costs $n$ queries.
\end{definition}
\begin{definition}\label{def:zr}
Write $g^t=\nabla f_{i_t}(x^t)$ and $h^t=\nabla^2f_{i_t}(x^t)v^t$ for the responses to the $t$-th query, and let $S^t=\bigcup_{s\le t}(\supp(g^s)\cup\supp(h^s))$. The first query is at $x^1=0$ with $v^1=0$. A zero-respecting algorithm satisfies $\supp(x^t)\cup\supp(v^t)\subseteq S^{t-1}$ for every $t\ge1$, and $\supp(\hat x^t)\subseteq S^t$ for every output $\hat x^t$ after $t$ queries.
\end{definition}
Query indices, points, directions, and stopping decisions may depend measurably on the transcript and private randomness. Linear-span methods started at the origin satisfy Definition~\ref{def:zr}. Methods using random starting vectors, including Oja's iteration and shift-and-invert routines, need not. We call algorithms without the support restriction \emph{unrestricted randomized} algorithms. For a function class $\mathcal F$, let $\N_\eps(\mathcal F)$ be the smallest fixed query budget with which some unrestricted randomized algorithm satisfies \eqref{eq:target} with probability at least $2/3$ on every instance of $\mathcal F$. We define $\Nzr_\eps(\mathcal F)$ in the same way for randomized zero-respecting algorithms, requiring only the gradient condition in \eqref{eq:target}, since a zero-respecting algorithm cannot certify curvature at a saddle point with zero gradient. Algorithms may know the class parameters and the dimension. Each of our lower bounds holds in a sufficiently large finite dimension and already prevents the gradient condition.

\subsection{Related Work}\label{sec:related}
\paragraph{Approximate local minima of finite sums.}
Perturbed gradient methods~\cite{jin} and negative-curvature extraction~\cite{neon2}, combined with nested or recursive variance reduction, find approximate local minima using $\tO(\sqrt n\,\eps^{-2}+n\eps^{-3/2}+n^{3/4}\eps^{-7/4})$ oracle calls in normalized parameters~\cite{snvrg,ssrgd}. For $\eps\le L^2/(n\rho)$, all of these bounds and \eqref{eq:fastcubic} are dominated by $n^{3/4}L^{1/2}\rho^{1/4}\Delta\eps^{-7/4}$ or by the larger term $\sqrt n\,L\Delta\eps^{-2}$. To our knowledge, previous lower bounds in this regime gave at most the sum of $\Omega(n)$ and the single-function rate.

\paragraph{Lower bounds for a single function.}
\citet{cdhs} established the lower bound $\Omega(\eps^{-12/7})$ for first-order methods on functions with Lipschitz Hessians. \citet{zhousharp} obtained the tight deterministic bound $\Omega(L^{1/2}\rho^{1/4}\Delta\eps^{-7/4})$ by a block chain, and \citet[Section~3.3]{cjz} gave an alternative proof by a nested chain whose hyper-objective is exactly the nonconvex chain of \citet{cdhs}. \citet{zhangrand} proved the same bound for randomized algorithms with values and gradients by a spectral construction. We embed the nested chain of \citet{cjz} into finite sums and extend the construction of \citet{zhangrand} to exact Hessian--vector products and replicate it across components.

\paragraph{Random edge owners and the factor $n^{3/4}$.}
For strongly convex sums of nonconvex components, the factor $n^{3/4}$ appears under averaged smoothness~\cite[Theorem~4.1]{zhougu} and is attained by Katyusha~X~\cite{katx}. \citet{han} obtain related bounds by partitioning tridiagonal chains among components. For nonconvex--strongly-concave finite sums, \citet{zhouncsc} assigns the internal edges of dual paths to random components and obtains the factor $n^{3/4}\sqrt\kappa$ per primal link. Our construction applies the same assignment to the quadratic bridges of a nested chain. Under individual smoothness, all of these constructions give only $\sqrt n$, and whether Katyusha~X is optimal in that setting is open.

\paragraph{Low accuracy and dense weak hiding.}
Without the curvature condition, recursive gradient estimators attain $O(n+\sqrt n\,L\Delta\eps^{-2})$ oracle calls under averaged smoothness~\cite{spider,page}. This rate is optimal under averaged smoothness~\cite{spider,zhougu} and, by the dense-weak-hiding construction of \citet{densehiding}, under individual smoothness. These results do not control the Hessian-Lipschitz constant at scale $\rho$. We show that the average of the dense-weak-hiding instances satisfies Assumption~\ref{ass:hess} when $\eps\gtrsim L^2/(n\rho)$ (Theorem~\ref{thm:low}).

\paragraph{Other oracle models.}
\citet{arjevani} study stochastic second-order oracles with fresh samples, whose rates do not transfer to finite sums. \citet{chensample} extend this model to reusable samples with full Hessians and leave the Hessian--vector model open. When each call returns a full component Hessian, cubic regularization~\cite{nesterov} needs only $O(1+\rho^{1/2}\Delta\eps^{-3/2})$ full Hessians, which changes the high-accuracy rate (Proposition~\ref{prop:hessian}). The matching lower bound uses components whose Hessian-Lipschitz constants grow with $n$, since Assumption~\ref{ass:hess} constrains only the average Hessian. When every component Hessian is $\rho$-Lipschitz, the recursive Taylor estimators of \citet{zhougu22} find a point with $\norm{\nabla f(\hat x)}\le\eps$ using $\tO(n+n^{3/4}\rho^{1/2}\Delta\eps^{-3/2})$ full component Hessians, and \citet{wu26} show that this rate is optimal for unrestricted randomized algorithms when no bound on $L$ is imposed; their characterization $\Theta_p(n+\Delta L_p^{1/p}n^{1-1/(2p)}\eps^{-(p+1)/p})$ covers oracles returning all derivatives through order $p\ge2$. Variance-reduced cubic Newton methods have also been analyzed under averaged Hessian assumptions~\cite{re3mcn}. Tables~\ref{tab:comparison} and~\ref{tab:low} compare the known bounds with ours.

\begin{table}[t]
\caption{The component-oracle complexity of achieving \eqref{eq:target} in the high-accuracy regime $\eps\le c_0L^2/(n\rho)$ with $\rho^{1/2}\Delta\eps^{-3/2}\ge C_0$. ZR denotes zero-respecting algorithms (Definition~\ref{def:zr}); all other rows allow unrestricted randomized algorithms with Hessian--vector products. $^\dagger$The single-function bounds apply to finite sums with identical components and do not grow with $n$. $^\ddagger$The IS lower bounds hold for $\eps\le c_0L^2/(n^2\rho)$. In this regime both FastCubic rows equal \eqref{eq:fastcubic} up to constants.}
\label{tab:comparison}
\centering
\footnotesize
\renewcommand{\arraystretch}{1.15}
\begin{tabular}{lcc}
\toprule
Methods & Averaged Smoothness & Individual Smoothness \\
\midrule
FastCubic \cite{fastcubic} & -- & $\tO\bigl(n^{3/4}L^{1/2}\rho^{1/4}\Delta\eps^{-7/4}\bigr)$ \\
FastCubic (Theorem~\ref{thm:fastcubic}) & $\tO\bigl(n^{3/4}L^{1/2}\rho^{1/4}\Delta\eps^{-7/4}\bigr)$ & $\tO\bigl(n^{3/4}L^{1/2}\rho^{1/4}\Delta\eps^{-7/4}\bigr)$ \\
\midrule
Lower Bounds, $n=1$ \cite{zhousharp,cjz,zhangrand}$^\dagger$ & $\Omega\bigl(L^{1/2}\rho^{1/4}\Delta\eps^{-7/4}\bigr)$ & $\Omega\bigl(L^{1/2}\rho^{1/4}\Delta\eps^{-7/4}\bigr)$ \\
Lower Bounds, ZR (Theorem~\ref{thm:nested}) & $\Omega\bigl(n^{3/4}L^{1/2}\rho^{1/4}\Delta\eps^{-7/4}\bigr)$ & $\Omega\bigl(\sqrt n\,L^{1/2}\rho^{1/4}\Delta\eps^{-7/4}\bigr)^\ddagger$ \\
Lower Bounds (Theorem~\ref{thm:highlower}) & $\Omega\bigl(n^{3/4}L^{1/2}\rho^{1/4}\Delta\eps^{-7/4}\bigr)$ & $\Omega\bigl(\sqrt n\,L^{1/2}\rho^{1/4}\Delta\eps^{-7/4}\bigr)^\ddagger$ \\
\bottomrule
\end{tabular}
\end{table}

\begin{table}[t]
\caption{The component-oracle complexity of achieving \eqref{eq:target} in the low-accuracy regime $L^2/(n\rho)\le\eps\le L^2/\rho$. AS and IS denote averaged and individual smoothness (Assumptions~\ref{ass:as} and~\ref{ass:is}). $^\dagger$These methods only guarantee $\norm{\nabla f(\hat x)}\le\eps$. $^\ddagger$Each oracle call returns the full component Hessian $\nabla^2f_i(x)$. $^\S$The bound assumes that every component Hessian is $\rho$-Lipschitz, a subclass of IS.}
\label{tab:low}
\centering
\footnotesize
\renewcommand{\arraystretch}{1.15}
\begin{tabular}{lcc}
\toprule
Methods & Assumption & Oracle Complexity \\
\midrule
SPIDER \cite{spider}, PAGE \cite{page}$^\dagger$ & AS & $O\bigl(n+\sqrt n\,L\Delta\eps^{-2}\bigr)$ \\
Cubic regularization \cite{nesterov}$^\ddagger$ & AS & $O\bigl(n+n\rho^{1/2}\Delta\eps^{-3/2}\bigr)$ \\
Recursive Taylor estimator \cite{zhougu22}$^{\dagger\ddagger\S}$ & IS & $\tO\bigl(n+n^{3/4}\rho^{1/2}\Delta\eps^{-3/2}\bigr)$ \\
FastCubic \cite{fastcubic} & IS & $\tO\bigl(n\rho^{1/2}\Delta\eps^{-3/2}+n^{3/4}L^{1/2}\rho^{1/4}\Delta\eps^{-7/4}\bigr)$ \\
PathVR (Algorithm~\ref{alg:pathvr}) & AS & $\tO\bigl(n+\sqrt n\,L\Delta\eps^{-2}\bigr)$ \\
\midrule
Lower Bound (Theorem~\ref{thm:low})$^\ddagger$ & IS & $\Omega\bigl(n+\sqrt n\,L\Delta\eps^{-2}\bigr)$ \\
\bottomrule
\end{tabular}
\end{table}

\section{Upper Bounds}\label{sec:upper}
We first develop a curvature search under averaged smoothness in Section~\ref{sec:curv}. We then introduce PathVR, which attains the low-accuracy upper bound, in Section~\ref{sec:pathvr}. Finally, Section~\ref{sec:fastcubic} extends FastCubic to averaged smoothness.

\subsection{Curvature Search under Averaged Smoothness}\label{sec:curv}
Given a point $x$ and a threshold $\delta$, we seek either a unit direction of curvature at most $-\delta/2$ or a certificate that $\lambda_{\min}(\nabla^2f(x))\ge-\delta$. Algorithm~\ref{alg:curv} applies Oja's iteration to the sampled shifted Hessians $2LI-\nabla^2f_{i_t}(x)$. Their mean $\Sigma=2LI-\nabla^2f(x)$ satisfies $LI\preceq\Sigma\preceq3LI$, and their variance is at most $L^2$ by \eqref{eq:asmatrix}. A top eigenvector of $\Sigma$ is therefore a bottom eigenvector of $\nabla^2f(x)$.

We validate each candidate using a median of sample means of its Rayleigh quotient, with fresh indices and the same variance bound. The algorithm returns the candidate with the smallest estimate if that estimate is sufficiently negative, and the zero vector otherwise. We choose the parameters
\begin{equation}\label{eq:curv-params}
 \begin{gathered}
 \alpha=\dfrac{\delta}{24L},\qquad
 \Lambda_d=\log\dfrac{ed}{\alpha},\qquad
 a=\dfrac{C_1\Lambda_d}{\alpha},\qquad
 T=\Bigl\lceil\dfrac{C_2\Lambda_d^2}{\alpha^2}\Bigr\rceil,\\
 R=\Bigl\lceil\dfrac{\log(2/p)}{\log4}\Bigr\rceil,\qquad
 m=\Bigl\lceil\dfrac{4096L^2}{\delta^2}\Bigr\rceil,\qquad
 M=\lceil C_3\log(2R/p)\rceil,
 \end{gathered}
\end{equation}
where $C_1,C_2,C_3$ are sufficiently large universal constants.

\begin{algorithm}[t]
\caption{$\mathrm{CurvSearch}(x;\delta,p)$}\label{alg:curv}
\begin{algorithmic}[1]
\Statex \textbf{Input:} point $x$, threshold $0<\delta\le L$, failure probability $0<p<1/2$; parameters \eqref{eq:curv-params}.
\For{$r=1,\ldots,R$}
 \State Draw $w_r^0$ uniformly from the unit sphere.
 \For{$t=0,1,\ldots,T-1$}
  \State Draw $i_t$ uniformly from $[n]$ and set $\eta_t\gets a/(3L(T/4+t))$.
  \State $\tilde w\gets w_r^t+\eta_t\bigl(2Lw_r^t-\nabla^2f_{i_t}(x)w_r^t\bigr)$.
  \If{$\norm{\tilde w}=0$}
   \State $w_r^T\gets e_1$; \textbf{break} the inner loop.
  \Else
   \State $w_r^{t+1}\gets\tilde w/\norm{\tilde w}$.
  \EndIf
 \EndFor
 \State $\hat\lambda_r\gets$ the median of $M$ independent averages, each of $m$ samples $(w_r^T)^\top\nabla^2f_i(x)w_r^T$ with $i$ uniform on $[n]$.
\EndFor
\State $r^\star\gets\arg\min_{r\in[R]}\hat\lambda_r$.
\State \textbf{if} $\hat\lambda_{r^\star}\le-5\delta/8$ \textbf{then return} $u=w_{r^\star}^T$ \textbf{else return} $u=0$.
\end{algorithmic}
\end{algorithm}

The gap-free analysis of Oja's iteration by \citet{gkty} requires only the variance bound $L^2$. It gives the following guarantee, proved in Appendix~\ref{app:nc}.
\begin{proposition}[cf. Gu et al.~{\cite[Lemma~10]{gkty}}]\label{prop:nc}
Assume \eqref{eq:asmatrix} at $x$, and let $0<\delta\le L$ and $0<p<1/2$. Algorithm~\ref{alg:curv} uses $O\bigl(\frac{L^2}{\delta^2}\bigl[\log^2\frac{edL}{\delta}+\log\frac{2\log(2/p)}{p}\bigr]\log\frac2p\bigr)$ component Hessian--vector products at $x$. With probability at least $1-p$, its output $u$ is either a unit vector with $u^\top\nabla^2f(x)u\le-\delta/2$, or $u=0$ and $\lambda_{\min}(\nabla^2f(x))\ge-\delta$.
\end{proposition}

\subsection{Path-Length-Based Variance Reduction}\label{sec:pathvr}
We build on the recursive gradient estimator of SPIDER~\cite{spider},
\begin{equation}\label{eq:spider}
 v^{t+1}=v^t+\dfrac1{b^t}\sum_{j\in\mathcal B^t}\bigl(\nabla f_j(x^{t+1})-\nabla f_j(x^t)\bigr),
\end{equation}
where $\mathcal B^t$ is a mini-batch of $b^t$ indices drawn uniformly from $[n]$. Under Assumption~\ref{ass:as}, the update \eqref{eq:spider} adds an error of variance at most $L^2\norm{x^{t+1}-x^t}^2/b^t$. The accuracy of $v^t$ thus depends on the \emph{lengths} of the steps since the last full gradient. SPIDER uses a fixed number of iterations between refreshes. A negative-curvature step has length $(\eps/\rho)^{1/2}$, which is much larger than the gradient-step length $\eps/(2L)$ when $\eps\le L^2/\rho$, and forces a full gradient at a cost of $n$ oracle calls.

We propose PathVR (path-length-based variance reduction, Algorithm~\ref{alg:pathvr}), which chooses both the mini-batch size and the refresh interval by step length. Given a window length $\ell>0$, a step of length $s^t=\norm{x^{t+1}-x^t}$ uses a mini-batch of size $b^t=\lceil ns^t/\ell\rceil$. We recompute the full gradient when the accumulated length $W^t$ of the current window would exceed $\ell$, or after a single step longer than $\ell/2$.

At each iteration, the method takes a normalized gradient step of length $\eta$ if $\norm{v^t}>2\eps$. Otherwise, it calls $\mathrm{CurvSearch}(x^t;\delta,p/(2K))$ from Algorithm~\ref{alg:curv}. If the search returns $u^t=0$, the method returns $x^t$, an approximate local minimum by Proposition~\ref{prop:nc}. Otherwise, it takes a step of length $r$ along the unit vector $u^t$, with its sign chosen against $v^t$. The tie-breaking rule $\sigma^t=1$ is necessary: without it, the method stalls at the saddle point of $f_i(x)=\cos x$ with $x^0=0$.

On the event of Theorem~\ref{thm:pathvr} below, the method terminates when the curvature search returns zero. The fallback output $x^0$ and the abort when the query budget $B$ is exhausted occur only on the failure event. They ensure a deterministic bound on the number of oracle calls.

\begin{algorithm}[t]
\caption{PathVR: Path-Length-Based Variance Reduction}\label{alg:pathvr}
\begin{algorithmic}[1]
\Statex \textbf{Input:} $\eps,\eta,\delta,r,\ell,K,B,p$; initial point $x^0=0$.
\State $v^0\gets\nabla f(x^0)$, $W^0\gets0$; abort with output $x^0$ once $B$ oracle calls have been made.
\For{$t=0,1,\ldots,K-1$}
 \If{$\norm{v^t}>2\eps$}
  \State $x^{t+1}\gets x^t-\eta\,v^t/\norm{v^t}$.
 \Else
  \State $u^t\gets\mathrm{CurvSearch}(x^t;\,\delta,\,p/(2K))$ by Algorithm~\ref{alg:curv}.
  \State \textbf{if} $u^t=0$ \textbf{then return} $x^t$.
  \State $\sigma^t\gets-\sign\langle v^t,u^t\rangle$, with $\sigma^t\gets1$ if $\langle v^t,u^t\rangle=0$.
  \State $x^{t+1}\gets x^t+\sigma^tru^t$.
 \EndIf
 \State $s^t\gets\norm{x^{t+1}-x^t}$, $b^t\gets\lceil ns^t/\ell\rceil$.
 \If{$s^t>\ell/2$ or $W^t+s^t>\ell$}
  \State $v^{t+1}\gets\nabla f(x^{t+1})$, $W^{t+1}\gets0$.
 \Else
  \State Draw $\mathcal B^t=\{j^t_1,\ldots,j^t_{b^t}\}$ uniformly and independently from $[n]$.
  \State $v^{t+1}\gets v^t+\frac1{b^t}\sum_{j\in\mathcal B^t}\bigl(\nabla f_j(x^{t+1})-\nabla f_j(x^t)\bigr)$, $W^{t+1}\gets W^t+s^t$.
 \EndIf
\EndFor
\State \Return $x^0$.
\end{algorithmic}
\end{algorithm}

The following lemma gives uniform accuracy of the estimator along the trajectory. We prove it using the Freedman inequality in Hilbert spaces~\cite[Theorem~2.1]{pzz}; under individual smoothness, the Azuma--Hoeffding inequality suffices.
\begin{lemma}\label{lem:pathvr-est}
Under Assumption~\ref{ass:as}, run Algorithm~\ref{alg:pathvr} with $p\in(0,1)$, $K\ge2$, $c_p=64\log^2(4K/p)$, and $\ell=\sqrt n\,\eps/(16\sqrt{c_p}\,L)$. Then with probability at least $1-p/2$, it holds that $\norm{v^t-\nabla f(x^t)}\le\eps/16$ for all iterations $t$ performed by the algorithm.
\end{lemma}
We choose the parameters
\begin{equation}\label{eq:pathvr-params}
 \begin{gathered}
 \eta=\dfrac{\eps}{2L},\qquad
 \delta=(\rho\eps)^{1/2},\qquad
 r=\Bigl(\dfrac{\eps}{\rho}\Bigr)^{1/2},\qquad
 \ell=\dfrac{\sqrt n\,\eps}{16\sqrt{c_p}\,L},\\
 K=\max\!\left\{2,\Bigl\lceil1+\dfrac{32L\Delta}{27\eps^2}+97\rho^{1/2}\Delta\eps^{-3/2}\Bigr\rceil\right\}.
 \end{gathered}
\end{equation}
with $c_p=64\log^2(4K/p)$, and set the query budget to
\begin{equation}\label{eq:pathvr-cap}
 B=\Bigl\lceil C\Bigl[n+\sqrt{c_pn}\,L\Delta\eps^{-2}+\bigl(1+\rho^{1/2}\Delta\eps^{-3/2}\bigr)Q_{\mathrm{nc}}\Bigr]\Bigr\rceil,
\end{equation}
where $C$ is a sufficiently large universal constant and $Q_{\mathrm{nc}}$ is the number of oracle calls in Proposition~\ref{prop:nc} with threshold $\delta$ and failure probability $p/(2K)$. On the event of Lemma~\ref{lem:pathvr-est} and successful curvature searches, each gradient step decreases $f$ by at least $\frac{27}{32}\eps^2/L$ and each curvature step by at least $\eps^{3/2}/(48\rho^{1/2})$. These descent bounds give the following theorem.
\begin{theorem}\label{thm:pathvr}
Under Assumptions~\ref{ass:as} and~\ref{ass:hess}, let $0<\eps\le L^2/\rho$ and $p\in(0,1)$, and run Algorithm~\ref{alg:pathvr} from $x^0=0$ with the parameters \eqref{eq:pathvr-params} and the budget \eqref{eq:pathvr-cap}. Then with probability at least $1-p$, the output $\hat x$ satisfies
\[
 \norm{\nabla f(\hat x)}\le\dfrac{33}{16}\eps
 \qquad\text{and}\qquad
 \lambda_{\min}\bigl(\nabla^2f(\hat x)\bigr)\ge-(\rho\eps)^{1/2},
\]
and the algorithm uses at most
\[
 \tO\!\left(n+\sqrt n\,L\Delta\eps^{-2}+\dfrac{L^2}{\rho\eps}+L^2\rho^{-1/2}\Delta\eps^{-5/2}\right)
\]
component oracle calls. The number of oracle calls never exceeds $B$, which is of the same order.
\end{theorem}
An amortized count gives the complexity bound. Both types of step decrease $f$ by at least $\eps s^t/48$, so the total path length is at most $48\Delta/\eps$. A step of length $s^t$ samples at most $1+ns^t/\ell$ component indices, each requiring two gradient evaluations. Each refresh costs $n$ calls and is triggered only after an accumulated length of at least $\ell/2$. The total cost of mini-batches and refreshes is therefore $O(n+K+n\Delta/(\eps\ell))=\tO(n+\sqrt n\,L\Delta\eps^{-2})$, giving the first two terms. The remaining two terms come from at most $1+48\rho^{1/2}\Delta\eps^{-3/2}$ curvature searches, each using $\tO(L^2/(\rho\eps))$ Hessian--vector products. In the low-accuracy regime $\eps\ge L^2/(n\rho)$, the third term is absorbed by $n$, and the fourth is at most a constant multiple of the second.
\begin{corollary}\label{cor:pathvr}
In the setting of Theorem~\ref{thm:pathvr} with $\eps=16\eps_{\mathrm{out}}/33$, where $L^2/(n\rho)\le\eps_{\mathrm{out}}\le L^2/\rho$, Algorithm~\ref{alg:pathvr} achieves \eqref{eq:target} with tolerance $\eps_{\mathrm{out}}$ and probability at least $1-p$ within a fixed budget of $\tO\bigl(n+\sqrt n\,L\Delta\eps_{\mathrm{out}}^{-2}\bigr)$ component oracle calls.
\end{corollary}
The proofs of Lemma~\ref{lem:pathvr-est}, Theorem~\ref{thm:pathvr}, and Corollary~\ref{cor:pathvr} are given in Appendix~\ref{app:pathvr}.

\subsection{FastCubic under Averaged Smoothness}\label{sec:fastcubic}
To extend FastCubic to averaged smoothness, we replace its inverse solver and eigenvector routine, the only two procedures that use individual smoothness. Both replacements require only the second-moment bound \eqref{eq:asmatrix}. We first give an inverse solver based on a variant of Katyusha~X whose analysis depends only on the variance.
\begin{lemma}[cf. Allen-Zhu~{\cite[Lemma~3.3]{katx}}]\label{lem:inverse}
Let $A=n^{-1}\sum_iH_i+aI\succeq\mu I$ with $\mu>0$, where $H_i$ satisfy \eqref{eq:asmatrix}. Let $\eta>0$, $p\in(0,1)$, and $\varkappa=(L+|a|)/\mu$. A randomized algorithm using $\tO(n+n^{3/4}\varkappa^{1/2})$ component matrix--vector products returns $z$ with $\norm{z-A^{-1}b}\le\eta\norm b$ with probability at least $1-p$.
\end{lemma}
By \eqref{eq:asmatrix}, the SVRG estimator has centered variance at most $L^2\norm{z-z_0}^2$, as required by the variance-only epoch lemma of \citet{katx}. The eigenvector routine applies AppxPCA~\cite[Theorem~3.1]{lazysvd} to the shifted inverses covered by Lemma~\ref{lem:inverse}. Since the outer loop of FastCubic uses only the average Hessian~\cite[Sections~3--7]{fastcubic}, these replacements give the following theorem. Its proof is in Appendix~\ref{app:fastcubic}.
\begin{theorem}\label{thm:fastcubic}
For $0<\eps\le L^2/\rho$ and $p\in(0,1)$, a variant of FastCubic achieves \eqref{eq:target} on $\Fas(n,L,\rho,\Delta)$ with probability at least $1-p$ within a fixed budget of
\begin{equation}\label{eq:fastcubic-as}
 \tO\!\left(\bigl(1+\rho^{1/2}\Delta\eps^{-3/2}\bigr)\bigl(n+n^{3/4}L^{1/2}\rho^{-1/4}\eps^{-1/4}\bigr)\right)
\end{equation}
queries. If $\rho^{1/2}\Delta\eps^{-3/2}\ge1$, this bound equals \eqref{eq:fastcubic}.
\end{theorem}
In the high-accuracy regime $\eps\le L^2/(n\rho)$, the second term of \eqref{eq:fastcubic} dominates. Section~\ref{sec:construction} gives a matching lower bound.

\section{Lower Bounds at High Accuracy}\label{sec:construction}
We next establish lower bounds in the high-accuracy regime. Section~\ref{sec:nestedsub} constructs a nested chain with random bridge owners for zero-respecting algorithms. Section~\ref{sec:highlower} extends the bound to unrestricted randomized algorithms through a replicated spectral construction.

\subsection{Nested Chain with Random Bridge Owners}\label{sec:nestedsub}
Let $\Upsilon(t)=120\int_1^t s^2(s-1)/(1+s^2)\,ds$ be the scalar potential of \citet{cdhs}. For integers $T,m\ge2$, let $\bar x\in\R^{T+1}$ and $x^{(i)}\in\R^m$ for $i\in[T]$. Following \citet[Equation~(10)]{cjz}, we define
\begin{equation}\label{eq:nested}
\begin{aligned}
 \hat f(x)&=\dfrac\nu2(1-\bar x_1)^2+\dfrac12\sum_{i=1}^{T}\hat Q_m\bigl(\bar x_i,x^{(i)},\bar x_{i+1}\bigr)+\nu\sum_{i=1}^{T}\Upsilon(\bar x_i),\\
 \hat Q_m(a,y,b)&=\dfrac1m(y_1-\theta a)^2+\sum_{j=1}^{m-1}(y_{j+1}-y_j)^2+\dfrac1m(y_m-\theta b)^2,
\end{aligned}
\end{equation}
with $m=\lfloor(3\nu)^{-1/2}\rfloor$ and $\theta=(\nu(3m-1))^{1/2}$. If $\nu\le1/12$, then $m\ge2$ and $\theta\sqrt m\le1$.

Each bridge $\hat Q_m$ connects consecutive outer coordinates and satisfies $\min_y\hat Q_m(a,y,b)=\nu(a-b)^2$ by \citet[Lemma~3.4]{cjz}. Thus, the \emph{reduced objective} obtained by minimizing over the bridge coordinates is exactly $\nu$ times the nonconvex chain of \citet{cdhs}. The full gradient norm is at least half the reduced-objective gradient norm~\cite[Lemma~3.5]{cjz}. In the coordinate order $\bar x_1,x^{(1)}_1,\ldots,x^{(1)}_m,\bar x_2,\ldots,\bar x_T,x^{(T)}_1,\ldots,x^{(T)}_m,\bar x_{T+1}$ given in \eqref{eq:order}, $\hat f$ is a zero-chain with gradient norm at least $\nu/8$ whenever the last two outer coordinates vanish. Since only the outer coordinates carry non-quadratic terms, extending the quadratic bridges does not increase their contribution to the Hessian-Lipschitz constant.

Write $d_{i,j}$ for the vector $e_k-e_{k'}$ that takes the difference of $x^{(i)}_j$ and $x^{(i)}_{j+1}$ for $j\in[m-1]$. For each bridge $i$, independently draw the owner $\sigma_{i,1}$ of its first edge uniformly from $[n]$. For $j\ge2$, draw $\sigma_{i,j}$ uniformly from $[n]\setminus\{\sigma_{i,j-1}\}$, so that equal-owner edges have disjoint supports. For scalings $\delta,\beta>0$, let $C(x)$ be $\delta\beta^2$ times the terms of $\hat f(x/\beta)$ other than the internal edges, and define
\begin{equation}\label{eq:nestedinstance}
 P_k=\sum_{(i,j):\,\sigma_{i,j}=k}d_{i,j}d_{i,j}^\top,\qquad
 f_k(x)=C(x)+\dfrac{n\delta}{2}\,x^\top P_kx\qquad(k\in[n]).
\end{equation}
For every assignment, $n^{-1}\sum_kf_k(x)=\delta\beta^2\hat f(x/\beta)$. The disjoint supports give $P_k^2=2P_k$, while $\sum_kP_k$ is a union of path Laplacians with norm at most $4$.

We specify the parameters in Appendix~\ref{app:nestedparams}. Taking $\delta\nu\asymp(\rho\eps)^{1/2}$ and $\beta\asymp\eps/(\delta\nu)$ makes the average $\rho$-Hessian-Lipschitz and keeps its gradient norm above $\eps$ until the chain is traversed. Taking $T\asymp\rho^{1/2}\Delta\eps^{-3/2}$ bounds the initial gap by $\Delta$. Averaged smoothness permits $\delta=L/(8\sqrt n)$, whereas individual smoothness permits only $\delta=L/(8n)$: the former follows from $\frac1n\sum_k\norm{n\delta P_kd}^2\le8n\delta^2\norm d^2$, and the latter because $\norm{\nabla^2f_k}$ contains the term $2n\delta$. The bridge-length condition $\nu\le1/12$ becomes $\eps\le5L^2/(768\,n\rho)$ and $\eps\le5L^2/(768\,n^2\rho)$, respectively.

The coordinates are exposed in the order \eqref{eq:order}. Since the endpoint terms are common to all components, any query can expose a bridge endpoint. An internal edge can be exposed only by querying its owner at a point or direction that is nonzero at the current frontier. Before the frontier reaches the edge, both endpoints vanish in every query. The transcript therefore gives no information about the owner, which remains conditionally uniform among $n$ or $n-1$ candidates. Discovering the owner requires order $n$ queries.

Reaching $\bar x_T$ requires traversing $(T-1)(m-1)$ internal edges. A conditional Hoeffding argument converts the discovery costs into a fixed-budget lower bound (Lemma~\ref{lem:search}). Multiplying the number of bridges, the bridge length, and the search cost gives
\begin{equation}\label{eq:product}
 \Omega(nTm)
 =\Omega\!\left(n\cdot\rho^{1/2}\Delta\eps^{-3/2}\cdot\nu^{-1/2}\right)
 =\begin{cases}
 \Omega\!\left(n^{3/4}L^{1/2}\rho^{1/4}\Delta\eps^{-7/4}\right)&\text{(AS)},\\
 \Omega\!\left(\sqrt n\,L^{1/2}\rho^{1/4}\Delta\eps^{-7/4}\right)&\text{(IS)}.
 \end{cases}
\end{equation}
For $n=1$, \eqref{eq:product} reduces to the bound of \citet{cjz}. The finite-sum structure has two effects. Discovering each edge costs order $n$ queries, and scaling its curvature to $L/\sqrt n$ or $L/n$ shortens the bridges by a factor of $n^{-1/4}$ or $n^{-1/2}$. Under averaged smoothness, the resulting factor $n\cdot n^{-1/4}=n^{3/4}$ matches \eqref{eq:fastcubic}.

We formalize this argument in the following theorem and give its proof in Appendix~\ref{app:nested}.
\begin{theorem}\label{thm:nested}
There are universal constants $c_0,C_0,c_1>0$ such that, for every $n\ge1$ and $\rho^{1/2}\Delta\eps^{-3/2}\ge C_0$,
\begin{align}
 \Nzr_\eps\bigl(\Fas(n,L,\rho,\Delta)\bigr)&\ge c_1\,n^{3/4}L^{1/2}\rho^{1/4}\Delta\eps^{-7/4},
 &&0<\eps\le c_0L^2/(n\rho),\label{eq:nested-as}\\
 \Nzr_\eps\bigl(\Fis(n,L,\rho,\Delta)\bigr)&\ge c_1\,\sqrt n\,L^{1/2}\rho^{1/4}\Delta\eps^{-7/4},
 &&0<\eps\le c_0L^2/(n^2\rho).\label{eq:nested-is}
\end{align}
The hard components share a common non-quadratic part, differ only in quadratic terms, and are individually $\rho$-Hessian-Lipschitz. The bounds allow adaptive component selection and exact Hessian--vector products.
\end{theorem}

\subsection{Replicated Spectra for Unrestricted Algorithms}\label{sec:highlower}
The nested chain gives a lower bound under the zero-respecting restriction. An unrestricted randomized algorithm can query points with small components on unrevealed coordinates, and an exact reply then reveals their directions. To obtain a lower bound without this restriction, we use the spectral construction of \citet{zhangrand}, whose hardness is governed by Krylov degree, and replicate its spectrum across components.

With the parameters of \citet{zhangrand} at smoothness level $L_\star$ (Appendix~\ref{app:replicated}), the construction has $N$ blocks and a diagonal matrix $\Lambda$ whose inverse cannot be approximated by polynomials of degree below $m$. These parameters satisfy $Nm\asymp L_\star^{1/2}\rho^{1/4}\Delta\eps^{-7/4}$ in the high-accuracy range. We fix $Nn$ public, mutually orthogonal owner subspaces, draw independent Haar rotations $U_{si}$ and uniform unit vectors $u_{si}$, and assign owner $i$ in block $s$ the matrix $A_{si}=U_{si}\Lambda U_{si}^\top$ and forcing vector $b_{si}\propto U_{si}u_{si}$. We set $f_i(z)=G(z)+\frac n2\sum_{s=1}^{N}z_{si}^\top A_{si}z_{si}$, where $G$ is the gated chain of \citet{zhangrand} built from the block averages $A_s,b_s$ in \eqref{eq:repinstance}.

The average $F$ is a block-diagonal instance of \citet{zhangrand}. Its gap, smoothness, and Hessian-Lipschitz bounds therefore hold for every realization, and an $\eps$-stationary point must open every gate. Since the owner Hessians have orthogonal supports, individual smoothness permits $L_\star=L/(2n)$, while averaged smoothness permits $L_\star=L/(2\sqrt n)$. The third derivatives of $f_i$ and $F$ coincide, so the components are individually $\rho$-Hessian-Lipschitz.

We couple the true oracle with an auxiliary oracle that discloses each block's rotations and forcing vectors on a fixed schedule and censors couplings beyond the current block. A conditional-Haar argument confines the latent query vectors of owner $(s,i)$ to a Krylov span. Its degree is bounded by the number of calls to that owner since disclosure (Lemma~\ref{lem:span}). A phase contains only $c_{\mathrm{delay}}nm$ calls, so at most $c_{\mathrm{delay}}n$ \emph{heavy} owners receive $m$ or more calls. Each heavy owner carries at most $2/n$ of the squared readout mass. For the light owners, the projection estimate of \citet[Lemma~7.2]{zhangrand} controls the information in their low-degree Krylov spans. Every readout consequently remains below its gate threshold until the phase ends. The censored and true replies coincide, giving a lower bound of $\Omega(nNm)$ calls (Appendix~\ref{app:highlower}).

Replicating the full spectrum in every owner is essential to this argument. If each eigenvalue belonged to a single component, one call per owner after disclosure of $b$ would reveal that owner's part of $A^{-1}b$. A phase would then cost $O(n+m)$ calls instead of $nm$.

The following theorem states the resulting lower bound. We give its proof in Appendix~\ref{app:highlower}.
\begin{theorem}\label{thm:highlower}
There are universal constants $c,C_0,c_2>0$ such that, if $\rho^{1/2}\Delta\eps^{-3/2}\ge C_0$, then
\begin{align}
 \N_\eps\bigl(\Fas(n,L,\rho,\Delta)\bigr)&\ge c_2\,n^{3/4}L^{1/2}\rho^{1/4}\Delta\eps^{-7/4},
 &&\eps\le c\,L^2/(n\rho),\label{eq:highlower-as}\\
 \N_\eps\bigl(\Fis(n,L,\rho,\Delta)\bigr)&\ge c_2\,\sqrt n\,L^{1/2}\rho^{1/4}\Delta\eps^{-7/4},
 &&\eps\le c\,L^2/(n^2\rho).\label{eq:highlower-is}
\end{align}
The hard components are individually $\rho$-Hessian-Lipschitz, and the bounds hold for exact Hessian--vector products and arbitrary query points.
\end{theorem}

\section{Lower Bound at Low Accuracy}\label{sec:lowlower}
For the lower bound in Theorem~\ref{thm:low}, we use the nonconvex dense-weak-hiding instance of \citet[Appendices~C.1 and~C.4]{densehiding} without modification (Appendix~\ref{app:dwh}). The instance has $T$ stages. Each stage hides a sign vector behind an encoding table with rows weakly biased toward that vector. The bias is $\varrho\asymp n^{-1/2}$, so recovering one stage requires order $n$ distinct rows. Consecutive stages are linked by the flat gate of \citet[Lemma~A.1]{densehiding}, whose derivatives vanish below its threshold. The scalings of \citet{densehiding} give individual smoothness $L$ and gradient norm at least $8\eps$ until the final stage is revealed.

We verify class membership and extend the information argument to our oracle. First, the row averages equal the hidden vectors exactly. Thus, the average of the unscaled components is independent of $\varrho$, although the component derivatives scale as $1/\varrho$. A universal bound on the third derivative of the average (Appendix~\ref{app:third}) gives
\[
 \Lip\Bigl(\nabla^2\dfrac1n\sum_if_i\Bigr)\le c_H\dfrac{L^2}{n\eps},\qquad c_H\le414.
\]
Replacing $L$ by $L_{\mathrm{eff}}=\min\{L,(n\rho\eps/c_H)^{1/2}\}$ places the instance in $\Fis(n,L,\rho,\Delta)$. The bound $\Omega(n+\sqrt n\,L_{\mathrm{eff}}\Delta\eps^{-2})$ of \citet[Theorem~3.1]{densehiding} then yields Theorem~\ref{thm:low}.

Second, the information argument of \citet[Appendix~B]{densehiding} requires only that any reply which does not complete the current stage depend on the history and the queried rows of that stage. Every derivative of an unopened link vanishes at such a point. The argument therefore applies to replies containing derivatives of any fixed order, including full Hessians (Lemma~\ref{lem:jets}).

Combining the lower bound with Corollary~\ref{cor:pathvr} gives the following complexity characterization.
\begin{theorem}\label{thm:low}
There are universal constants $n_0,c_3>0$ such that, for $n\ge n_0$, $L^2/(n\rho)\le\eps\le L^2/\rho$, and $\eps^2\le c_3L\Delta$,
\begin{equation}\label{eq:lowcomplexity}
 \N_\eps\bigl(\Fis(n,L,\rho,\Delta)\bigr)
 =\N_\eps\bigl(\Fas(n,L,\rho,\Delta)\bigr)
 =\tTheta\bigl(n+\sqrt n\,L\Delta\eps^{-2}\bigr).
\end{equation}
The upper bound is attained by PathVR (Corollary~\ref{cor:pathvr}). The lower bound holds even when each call returns the full component Hessian, and more generally equals $\Omega(n+\min\{\sqrt n\,L\Delta\eps^{-2},n\rho^{1/2}\Delta\eps^{-3/2}\})$ (Proposition~\ref{prop:lowgeneral}).
\end{theorem}
For the terms proportional to $\Delta$, the FastCubic bound exceeds the PathVR bound by a factor of order $(n\rho\eps)^{1/2}/L$, up to logarithmic factors.

\section{Near-Optimal Complexity}\label{sec:results}
We now combine the upper and lower bounds. Theorems~\ref{thm:fastcubic}, \ref{thm:nested}, and~\ref{thm:highlower} give the following characterization of the high-accuracy regime. The proof is in Appendix~\ref{app:matching}.
\begin{theorem}\label{thm:matching}
There are universal constants $c_0,C_0,c_1>0$ such that, for every $n\ge1$, $0<\eps\le c_0L^2/(n\rho)$, and $\rho^{1/2}\Delta\eps^{-3/2}\ge C_0$,
\begin{equation}\label{eq:matchedcomplexity}
 \N_\eps\bigl(\Fas(n,L,\rho,\Delta)\bigr)=\tTheta\bigl(n^{3/4}L^{1/2}\rho^{1/4}\Delta\eps^{-7/4}\bigr),
 \qquad
 \Nzr_\eps\bigl(\Fas(n,L,\rho,\Delta)\bigr)\ge c_1n^{3/4}L^{1/2}\rho^{1/4}\Delta\eps^{-7/4}.
\end{equation}
The lower bounds in \eqref{eq:matchedcomplexity} carry no logarithmic factors.
\end{theorem}
The parameter range in Theorem~\ref{thm:matching} ensures that the quadratic bridges of Section~\ref{sec:nestedsub} contain at least two vertices. In the range covered by the individual-smoothness lower bound, the bounds differ by a factor of $n^{1/4}$.
\begin{theorem}\label{thm:is}
Under the conditions of Theorem~\ref{thm:highlower}, for $\eps\le c\,L^2/(n^2\rho)$,
\begin{equation}\label{eq:isgap}
 c_2\sqrt n\,L^{1/2}\rho^{1/4}\Delta\eps^{-7/4}
 \le\N_\eps\bigl(\Fis(n,L,\rho,\Delta)\bigr)
 \le\tO\bigl(n^{3/4}L^{1/2}\rho^{1/4}\Delta\eps^{-7/4}\bigr).
\end{equation}
\end{theorem}
The gap in \eqref{eq:isgap} is also present for linear systems with individually bounded, possibly indefinite component matrices. Katyusha~X uses $n+n^{3/4}\sqrt\kappa$ queries, while the best lower bound under individual smoothness is of order $n+\sqrt{n\kappa}$~\cite{zhougu,han}. A solver using $n+n^{3/4-\theta}\sqrt\kappa$ queries would improve the second term of \eqref{eq:fastcubic} to $n^{3/4-\theta}L^{1/2}\rho^{1/4}\Delta\eps^{-7/4}$.

We also consider two related oracle models (Appendices~\ref{app:matching} and~\ref{app:single}). If each call returns the full component Hessian, cubic regularization requires only $O(1+\rho^{1/2}\Delta\eps^{-3/2})$ full Hessians. The complexity is determined by the better of the cubic-regularization and PathVR bounds.
\begin{proposition}\label{prop:hessian}
Suppose each call may return the full matrix $\nabla^2f_i(x)$. For $n\ge n_0$ and in either accuracy regime under the corresponding conditions of Theorem~\ref{thm:matching} or~\ref{thm:low}, the fixed-budget complexity over $\Fas$ or $\Fis$ is $\tTheta(n+\min\{\sqrt n\,L\Delta\eps^{-2},\,n\rho^{1/2}\Delta\eps^{-3/2}\})$.
\end{proposition}
The term $n\rho^{1/2}\Delta\eps^{-3/2}$ reflects that Assumption~\ref{ass:hess} constrains only the average Hessian: the dense-weak-hiding components behind Proposition~\ref{prop:lowgeneral} have Hessian-Lipschitz constants of order $\sqrt n$ times that of their average (Appendix~\ref{app:dwh}). On the subclass of $\Fis$ whose components are individually $\rho$-Hessian-Lipschitz, \citet{zhougu22} find a point with $\norm{\nabla f(\hat x)}\le\eps$ using $\tO(n+n^{3/4}\rho^{1/2}\Delta\eps^{-3/2})$ full component Hessians, and \citet{wu26} prove the lower bound $\Omega(n+n^{3/4}\rho^{1/2}\Delta\eps^{-3/2})$ for unrestricted randomized algorithms over classes with no bound on $L$. Our high-accuracy hard instances are individually $\rho$-Hessian-Lipschitz (Theorems~\ref{thm:nested} and~\ref{thm:highlower}), so on this subclass the Hessian--vector complexity $n^{3/4}L^{1/2}\rho^{1/4}\Delta\eps^{-7/4}$ exceeds the full-Hessian rate by the factor $(L^2/(\rho\eps))^{1/4}\ge1$.
For a single function, the replicated spectral construction with $n=1$ extends the lower bound of \citet{zhangrand} to exact Hessian--vector products.
\begin{proposition}\label{prop:single}
For $n=1$, $\eps\le cL^2/\rho$, and $\rho^{1/2}\Delta\eps^{-3/2}\ge C_0$, every randomized algorithm with values, gradients, and exact Hessian--vector products needs $\Omega(L^{1/2}\rho^{1/4}\Delta\eps^{-7/4})$ queries, which matches FastCubic up to logarithmic factors.
\end{proposition}

\section{Conclusion}\label{sec:conclusion}
We characterize the component-oracle complexity of finding approximate local minima of nonconvex finite sums under averaged smoothness and a Lipschitz average Hessian. At high accuracy, the minimax complexity is $\tTheta(n^{3/4}L^{1/2}\rho^{1/4}\Delta\eps^{-7/4})$, and FastCubic attains the upper bound. The nested-chain construction explains the factor $n^{3/4}$ through the cost of discovering randomly owned bridge edges. At low accuracy, PathVR attains $\tO(n+\sqrt n\,L\Delta\eps^{-2})$, matching the dense-weak-hiding lower bound. In the range covered by the individual-smoothness lower bound, the bounds differ by a factor of $n^{1/4}$. Faster linear-system solvers that exploit individual smoothness would narrow this gap. Other open questions are to remove logarithmic factors from the upper bound, characterize the transition $\eps\asymp L^2/(n\rho)$ up to constants, and extend random-owner constructions to finite-sum minimax and bilevel problems with Lipschitz Hessians.

\section*{Acknowledgements}
The author thanks Lesi Chen of Tsinghua University for helpful discussions of this open problem.

Codex and Claude Science were used to assist with the development of the constructions and proofs, numerical checks of algebraic identities, the review of related literature, and the preparation of this manuscript. The author is responsible for the mathematical claims, references, and final presentation.

\begingroup
\fontsize{10}{12}\selectfont
\bibliographystyle{plainnat}
\bibliography{refs}

\endgroup

\clearpage
\appendix
\section{Proofs for Section~\ref{sec:upper}}\label{app:upper}
\subsection{Proof of Proposition~\ref{prop:nc}}\label{app:nc}
\begin{proof}
Write $H_i=\nabla^2f_i(x)$, $H=\nabla^2f(x)$, $A_t=2LI-H_{i_t}$, and $\Sigma=2LI-H$. By \eqref{eq:asmatrix},
\[
 LI\preceq\Sigma\preceq3LI,\qquad
 \E(A_t-\Sigma)^2\preceq L^2I,\qquad
 \dfrac{a\lambda_1(\Sigma)}{3L}\in[a/3,a].
\]
Thus \citet[Lemma~10]{gkty} applies to the step sizes in Algorithm~\ref{alg:curv} for sufficiently large $C_1,C_2$. Its good event has a nonzero final unnormalized product, and a zero intermediate update forces a zero final product, so a run assigned $e_1$ lies outside the good event; its validation is unchanged, and the success probability is unaffected. Each run produces a unit vector $w_r^T$ with
\[
 (w_r^T)^\top H w_r^T\le\lambda_{\min}(H)+\dfrac\delta8
\]
with probability at least $3/4$. Independence gives a successful run with probability at least $1-4^{-R}\ge1-p/2$.

Conditional on each candidate, fresh validation samples have variance at most $L^2$. Chebyshev's inequality gives
\[
 \Prob\left(\left|\dfrac1m\sum_{j=1}^m(w_r^T)^\top H_{i_j}w_r^T
 -(w_r^T)^\top H w_r^T\right|>\dfrac\delta{16}\right)
 \le\dfrac{256L^2}{m\delta^2}\le\dfrac1{16}.
\]
The median and a union bound make all $R$ estimates accurate to $\delta/16$ with probability at least $1-p/2$, for sufficiently large $C_3$. On the intersection, if $\hat\lambda_{r^\star}\le-5\delta/8$, then
\[
 u^\top H u\le\hat\lambda_{r^\star}+\dfrac\delta{16}
 \le-\dfrac{9\delta}{16}\le-\dfrac\delta2.
\]
Otherwise $u=0$, and a successful run $w_r^T$ gives
\[
 \lambda_{\min}(H)
 \ge(w_r^T)^\top H w_r^T-\dfrac\delta8
 \ge\hat\lambda_r-\dfrac{3\delta}{16}
 \ge\hat\lambda_{r^\star}-\dfrac{3\delta}{16}
 >-\dfrac{13\delta}{16}>-\delta.
\]
The cost is $RT+RMm$ Hessian--vector products, yielding the stated bound.
\end{proof}

\subsection{Path-Length-Based Variance Reduction}\label{app:pathvr}
Use the setting of Theorem~\ref{thm:pathvr}, with $s^t=\norm{x^{t+1}-x^t}$ and $b^t=\lceil ns^t/\ell\rceil$. Gradient steps, curvature steps, and refreshes refer to the normalized gradient update, the signed curvature update, and the full-gradient recomputation in Algorithm~\ref{alg:pathvr}. A \emph{window} consists of consecutive iterations between refreshes, including the initial window.

\subsubsection{Basic Lemmas}
\begin{lemma}[Peng et al.~{\cite[Theorem~2.1]{pzz}}]\label{lem:freedman}
Let $\{\xi_k\}_{k\ge1}$ be a martingale difference sequence in a Hilbert space with respect to a filtration $\{\mathcal F_k\}_{k\ge0}$ such that $\norm{\xi_k}\le M$ almost surely for all $k$. Let $V_m=\sum_{k=1}^m\E[\norm{\xi_k}^2\mid\mathcal F_{k-1}]$. Then for any $a>0$ and $V>0$, we have
\[
 \Prob\Bigl(\exists\,m\ge1:\ \Bigl\|\sum_{k=1}^m\xi_k\Bigr\|\ge a\ \text{and}\ V_m\le V\Bigr)
 \le2\exp\Bigl(-\dfrac{a^2}{2(V+aM/3)}\Bigr).
\]
\end{lemma}

\begin{lemma}\label{lem:as-facts}
Suppose Assumption~\ref{ass:as} holds. For any $x,x'\in\R^d$ with $s=\norm{x'-x}$, let $d_i=\nabla f_i(x')-\nabla f_i(x)$ for $i\in[n]$ and $\bar d=\nabla f(x')-\nabla f(x)$. Then we have
\begin{enumerate}
\item[(a)] $\norm{\bar d}\le Ls$ and $\max_{i\in[n]}\norm{d_i}\le\sqrt n\,Ls$;
\item[(b)] $\frac1n\sum_{i=1}^n\norm{d_i-\bar d}^2\le L^2s^2$;
\item[(c)] $f(x')\le f(x)+\langle\nabla f(x),x'-x\rangle+\frac L2s^2$.
\end{enumerate}
\end{lemma}
\begin{proof}
Jensen's inequality and \eqref{eq:as} give (a)--(b):
\[
 \begin{gathered}
 \norm{\bar d}^2 \le\dfrac1n\sum_i\norm{d_i}^2\le L^2s^2,\qquad \norm{d_i}^2 \le\sum_{i'}\norm{d_{i'}}^2\le nL^2s^2,\\
 \dfrac1n\sum_i\norm{d_i-\bar d}^2 =\dfrac1n\sum_i\norm{d_i}^2-\norm{\bar d}^2\le L^2s^2.
 \end{gathered}
\]
Part~(c) follows from $\norm{\nabla^2f}\le L$ in \eqref{eq:asmatrix}.
\end{proof}

\subsubsection{Proof of Lemma~\ref{lem:pathvr-est}}
\begin{lemma}\label{lem:window}
Under the setting of Lemma~\ref{lem:pathvr-est}, fix a window and let $t_0$ be its first iteration, so that $v^{t_0}=\nabla f(x^{t_0})$. For every iteration $t$ of the window that is not a refresh and every $k\in[b^t]$, define
\[
 \xi^t_k=\dfrac1{b^t}\Bigl(d^t_{j^t_k}-\bar d^{\,t}\Bigr),\qquad
 d^t_i=\nabla f_i(x^{t+1})-\nabla f_i(x^t),\qquad
 \bar d^{\,t}=\nabla f(x^{t+1})-\nabla f(x^t).
\]
Then, ordering the increments $\xi^t_k$ lexicographically in $(t,k)$, we have
\begin{enumerate}
\item[(a)] $v^{t+1}-\nabla f(x^{t+1})=\sum_{t'\le t}\sum_{k=1}^{b^{t'}}\xi^{t'}_k$ for every non-refresh iteration $t$ of the window;
\item[(b)] $\{\xi^t_k\}$ is a martingale difference sequence with $\norm{\xi^t_k}\le\eps/(8\sqrt{c_p})$ almost surely;
\item[(c)] the predictable quadratic variation of the whole window is at most $\eps^2/(256c_p)$.
\end{enumerate}
\end{lemma}
\begin{proof}
For each non-refresh iteration, the recursive estimator update in Algorithm~\ref{alg:pathvr} gives
\[
 v^{t+1}-\nabla f(x^{t+1}) =v^t-\nabla f(x^t)+\sum_{k=1}^{b^t}\xi^t_k =\sum_{t'=t_0}^{t}\sum_{k=1}^{b^{t'}}\xi^{t'}_k,
\]
since $v^{t_0}=\nabla f(x^{t_0})$. Let $\mathcal F^t_k$ be the history before sampling $j^t_k$. The points $x^t,x^{t+1}$ are already fixed, and $j^t_k$ is conditionally uniform and independent of earlier draws. Thus
\[
 \E[\xi^t_k\mid\mathcal F^t_k]
 =\dfrac1{b^t}\left(\dfrac1n\sum_id^t_i-\bar d^{\,t}\right)=0.
\]
Lemma~\ref{lem:as-facts} and $b^t\ge ns^t/\ell$ give
\[
 \norm{\xi^t_k} \le\dfrac{(\sqrt n+1)Ls^t}{b^t} \le\dfrac{2L\ell}{\sqrt n}=\dfrac\eps{8\sqrt{c_p}},\qquad \E[\norm{\xi^t_k}^2\mid\mathcal F^t_k] \le\dfrac{L^2(s^t)^2}{(b^t)^2}.
\]
The refresh rule ensures $\sum_ts^t\le\ell$ over the non-refresh iterations of a window; hence
\[
 \sum_t\sum_{k=1}^{b^t}\E[\norm{\xi^t_k}^2\mid\mathcal F^t_k] \le\sum_t\dfrac{L^2(s^t)^2}{b^t} \le\dfrac{L^2\ell}{n}\sum_ts^t \le\dfrac{L^2\ell^2}{n}=\dfrac{\eps^2}{256c_p}.
\]
\end{proof}

\begin{proof}[Proof of Lemma~\ref{lem:pathvr-est}]
Condition on a window's starting history and apply Lemmas~\ref{lem:freedman} and~\ref{lem:window} with
\[
 a=\dfrac\eps{16},\qquad M=\dfrac\eps{8\sqrt{c_p}},\qquad
 V=\dfrac{\eps^2}{256c_p},\qquad q=\log\dfrac{4K}{p}.
\]
Since $c_p=64q^2$ and $q\ge\log8>1$,
\[
 \dfrac{a^2}{2(V+aM/3)}
 =\dfrac1{2/c_p+4/(3\sqrt{c_p})}
 =\dfrac{192q^2}{32q+6}\ge5q.
\]
The probability of an error greater than $\eps/16$ anywhere in the window is therefore at most
\[
 2e^{-5q}=2\left(\dfrac p{4K}\right)^5\le\dfrac p{2K}.
\]
The estimator is exact at initialization and refreshes. A union bound over at most $K$ windows proves the claim.
\end{proof}

\subsubsection{Proof of Theorem~\ref{thm:pathvr}}
Let $\mathcal E$ be the event that every estimator error is at most $\eps/16$ and every call to $\mathrm{CurvSearch}$ is correct. Lemma~\ref{lem:pathvr-est} and a union bound over at most $K$ searches of failure probability $p/(2K)$ give $\Prob(\mathcal E)\ge1-p$.
\begin{lemma}\label{lem:descent}
Under the setting of Theorem~\ref{thm:pathvr}, the following statements hold on the event $\mathcal E$:
\begin{enumerate}
\item[(a)] if iteration $t$ is a gradient step, then $f(x^{t+1})\le f(x^t)-\dfrac{27\eps^2}{32L}$;
\item[(b)] if iteration $t$ is a curvature step, then $f(x^{t+1})\le f(x^t)-\dfrac{\eps^{3/2}}{48\rho^{1/2}}$;
\item[(c)] if the algorithm returns $x^t$ after the curvature search returns zero, then $\norm{\nabla f(x^t)}\le\frac{33}{16}\eps$ and $\lambda_{\min}(\nabla^2f(x^t))\ge-(\rho\eps)^{1/2}$.
\end{enumerate}
\end{lemma}
\begin{proof}
For a gradient step, Lemma~\ref{lem:as-facts}(c) gives
\[
 \begin{aligned}
 \left\langle\nabla f(x^t),\dfrac{v^t}{\norm{v^t}}\right\rangle
 &\ge\norm{v^t}-\norm{v^t-\nabla f(x^t)}
 \ge2\eps-\dfrac\eps{16}=\dfrac{31\eps}{16},\\
 f(x^{t+1})-f(x^t)
 &\le-\eta\dfrac{31\eps}{16}+\dfrac{L\eta^2}{2}
 =-\dfrac{27\eps^2}{32L}.
 \end{aligned}
\]
For a curvature step, $(u^t)^\top\nabla^2f(x^t)u^t\le-\delta/2$ and the sign rule yields
\[
 \sigma^t\langle\nabla f(x^t),u^t\rangle
 \le-\abs{\langle v^t,u^t\rangle}+\norm{\nabla f(x^t)-v^t}
 \le\dfrac\eps{16}.
\]
Taylor's bound from Assumption~\ref{ass:hess}, with $r=\delta/\rho$ and $\eps=\delta^2/\rho$, gives
\[
 f(x^{t+1})-f(x^t) \le\dfrac{r\eps}{16}-\dfrac{r^2\delta}{4}+\dfrac{\rho r^3}{6} =\dfrac{\delta^3}{\rho^2}\left(\dfrac1{16}-\dfrac14+\dfrac16\right) =-\dfrac{\eps^{3/2}}{48\rho^{1/2}}.
\]
At a return after the curvature search returns zero,
\[
 \norm{\nabla f(x^t)}\le\norm{v^t}+\dfrac\eps{16}\le\dfrac{33\eps}{16},
 \qquad \lambda_{\min}(\nabla^2f(x^t))\ge-\delta=-(\rho\eps)^{1/2}.
\]
\end{proof}

\begin{proof}[Proof of Theorem~\ref{thm:pathvr}]
On $\mathcal E$, temporarily ignoring the query cap, Lemma~\ref{lem:descent} and the gap bound give
\[
 N_{\mathrm g}\le\dfrac{32L\Delta}{27\eps^2},\qquad
 N_{\mathrm c}\le48\rho^{1/2}\Delta\eps^{-3/2},\qquad
 N_{\mathrm g}+N_{\mathrm c}+1\le K.
\]
Thus the algorithm returns after the curvature search returns zero with the required guarantees. Both types of step decrease the objective by at least $\eps s^t/48$, so its path length satisfies
\[
 P:=\sum_ts^t
 \le\dfrac{48}{\eps}\sum_t\bigl(f(x^t)-f(x^{t+1})\bigr)
 \le\dfrac{48\Delta}{\eps},
\]
and consequently
\begin{equation}\label{eq:nPl}
 \dfrac{nP}{\ell}
 \le\dfrac{48\Delta}{\eps}\cdot\dfrac{16\sqrt{c_p}\,L\sqrt n}{\eps}
 =768\sqrt{c_pn}\,L\Delta\eps^{-2}.
\end{equation}
Non-refresh iterations cost $2b^t\le2+2ns^t/\ell$ calls each, totaling at most $2K+2nP/\ell$. Initialization costs $n$. Refreshes with $s^t>\ell/2$ cost at most $2nP/\ell$ in total. Every other refresh closes a window with $W^t>\ell/2$; at most $2P/\ell$ such windows cost another $2nP/\ell$.

There are at most $1+N_{\mathrm c}$ curvature searches, each costing $Q_{\mathrm{nc}}=\tO(L^2/(\rho\eps))$ calls. By \eqref{eq:nPl}, the total is at most
\[
 \begin{aligned}
 &n+2K+\dfrac{6nP}{\ell}
   +\bigl(1+48\rho^{1/2}\Delta\eps^{-3/2}\bigr)Q_{\mathrm{nc}}\\
 &\quad\le n+2K+4608\sqrt{c_pn}\,L\Delta\eps^{-2}
   +\bigl(1+48\rho^{1/2}\Delta\eps^{-3/2}\bigr)Q_{\mathrm{nc}}.
 \end{aligned}
\]
Since $\eps\le L^2/\rho$,
\[
 \rho^{1/2}\Delta\eps^{-3/2}\le L\Delta\eps^{-2},\qquad
 K\le2+99L\Delta\eps^{-2}\le2+99\sqrt{c_pn}\,L\Delta\eps^{-2}.
\]
A sufficiently large $C$ in \eqref{eq:pathvr-cap} therefore makes the total at most
\[
 B=\tO\!\left(n+\sqrt n\,L\Delta\eps^{-2}
 +\dfrac{L^2}{\rho\eps}+L^2\rho^{-1/2}\Delta\eps^{-5/2}\right).
\]
Thus the cap is not reached on $\mathcal E$; outside $\mathcal E$, it enforces the same fixed budget.
\end{proof}

\subsubsection{Proof of Corollary~\ref{cor:pathvr}}
\begin{proof}
Set $\eps=16\eps_{\mathrm{out}}/33$. Theorem~\ref{thm:pathvr} gives, with probability at least $1-p$,
\[
 \norm{\nabla f(\hat x)}\le\dfrac{33\eps}{16}=\eps_{\mathrm{out}},\qquad
 \lambda_{\min}(\nabla^2f(\hat x))\ge-(\rho\eps)^{1/2}
 \ge-(\rho\eps_{\mathrm{out}})^{1/2}.
\]
Since $\eps_{\mathrm{out}}\ge L^2/(n\rho)$,
\[
 \dfrac{L^2}{\rho\eps} \le\dfrac{33n}{16},\qquad L^2\rho^{-1/2}\Delta\eps^{-5/2} =\left(\dfrac{L^2}{\rho\eps}\right)^{1/2}L\Delta\eps^{-2} \le\left(\dfrac{33}{16}\right)^{1/2}\sqrt n\,L\Delta\eps^{-2}.
\]
Substitution into Theorem~\ref{thm:pathvr} gives the budget $\tO(n+\sqrt n\,L\Delta\eps_{\mathrm{out}}^{-2})$.
\end{proof}

\subsection{Proof of Lemma~\ref{lem:inverse} and Theorem~\ref{thm:fastcubic}}\label{app:fastcubic}
\begin{proof}[Proof of Lemma~\ref{lem:inverse}]
Apply Katyusha~X to
\[
 q(z)=\dfrac12z^\top Az-b^\top z,\qquad
 q_i(z)=\dfrac12z^\top(H_i+aI)z-b^\top z.
\]
The snapshot estimator is unbiased and, by \eqref{eq:asmatrix},
\[
 \E_i\norm{(H_i-H)(z-z_0)}^2\le L^2\norm{z-z_0}^2.
\]
This verifies Lemma~3.3 of the full version of \citet{katx}; its acceleration argument uses convexity and smoothness only of the average quadratic. Set
\[
 S=L+|a|,\qquad m=\max\{n,2\},\qquad
 h=\min\left\{\dfrac1{2S},\dfrac1{2L\sqrt m}\right\}.
\]
By Lemma~4.1 and Theorem~4.3 of the full version of \citet{katx}, for $z_\star=A^{-1}b$ and initialization at zero,
\[
 \begin{gathered}
 \tau^{-1} =O(1+n^{-1/4}\varkappa^{1/2}),\qquad q(0)-q(z_\star) \le\dfrac{\norm b^2}{2\mu},\\
 \E[q(z_K)-q(z_\star)] \le2(1+\tau)^{-K}[q(0)-q(z_\star)].
 \end{gathered}
\]
Return zero if $b=0$ or $\mu\eta\ge1$. Otherwise take
\[
 K=\left\lceil\dfrac{\log(4/(p\mu^2\eta^2))}{\log(1+\tau)}\right\rceil.
\]
Strong convexity and Markov's inequality give
\[
 \Prob\bigl(\norm{z_K-z_\star}>\eta\norm b\bigr)
 \le\dfrac{2\E[q(z_K)-q(z_\star)]}{\mu\eta^2\norm b^2}
 \le\dfrac p2.
\]
Abort when a geometric epoch exceeds $\ell=\lceil m\log(2K/p)\rceil$ updates. The geometric-tail bound and a union bound give abort probability at most $p/2$. Coupling to the ideal execution therefore gives success probability at least $1-p$ within
\[
 K(n+2\ell)=\tO(n+n^{3/4}\varkappa^{1/2})
\]
queries.
\end{proof}

\begin{proof}[Proof of Theorem~\ref{thm:fastcubic}]
The shifts in FastCubicMin have averaged condition number $O(L(\rho\eps)^{-1/2})$ \cite[Appendix~G]{fastcubic}, so Lemma~\ref{lem:inverse} applies to its inverse--vector products. For the bottom eigenvector, apply AppxPCA \cite[Theorem~3.1]{lazysvd} with multiplicative and leakage tolerances $\xi/(8L)$, where
\[
 \xi=\dfrac{(\rho\eps)^{1/2}}{300},\qquad
 P=\dfrac{3LI-H}{4L},\qquad
 \dfrac I2\preceq P\preceq I.
\]
Its shifted inverses satisfy
\[
 \begin{gathered}
 \sigma I-P=\dfrac{H+(4L\sigma-3L)I}{4L},\qquad \lambda_{\min}\bigl(H+(4L\sigma-3L)I\bigr)\ge\dfrac\xi{48},\\
 |4L\sigma-3L|\le3L,\qquad \lambda_{\min}(\sigma I-P)\ge\dfrac\xi{192L}.
 \end{gathered}
\]
Thus Lemma~\ref{lem:inverse} applies without positive semidefinite component matrices. The outer loop uses only the average Hessian, its Lipschitz constant, and the prescribed inverse accuracies \cite[Sections~3--7]{fastcubic}. It solves
\[
 \begin{aligned}
 &O\!\left(1+\rho^{1/2}\Delta\eps^{-3/2}\right)
 \quad\text{cubic subproblems},\\
 &\tO\!\left(n+n^{3/4}L^{1/2}\rho^{-1/4}\eps^{-1/4}\right)
 \quad\text{queries per subproblem}.
 \end{aligned}
\]
Cap the outer steps at a sufficiently large universal multiple of this bound, distribute the failure probability across all subroutine calls, and abort before exceeding the total-query cap. Simultaneous subroutine success guarantees termination within both caps. Along successful descent histories, $\norm{\nabla f(x_t)}\le\sqrt{2L\Delta}$ controls the auxiliary logarithms. Multiplying the two bounds proves \eqref{eq:fastcubic-as}.
\end{proof}

\section{Proof of Theorem~\ref{thm:nested}}\label{app:nested}
\subsection{Parameters}\label{app:nestedparams}
Use the instance \eqref{eq:nested}--\eqref{eq:nestedinstance} with coordinate order
\begin{equation}\label{eq:order}
 \bar x_1,x^{(1)}_1,\ldots,x^{(1)}_m,\bar x_2,\ldots,
 \bar x_T,x^{(T)}_1,\ldots,x^{(T)}_m,\bar x_{T+1},
\end{equation}
and parameters
\begin{equation}\label{eq:nestedparams}
 \delta\nu=\Bigl(\dfrac{\rho\eps}{60}\Bigr)^{1/2},\qquad \beta=\dfrac{8\eps}{\delta\nu},\qquad T=\Bigl\lfloor\dfrac{\Delta\delta\nu}{1280\eps^2}\Bigr\rfloor,\qquad \delta=\begin{cases}L/(8\sqrt n)&\text{(AS)},\\ L/(8n)&\text{(IS)}.\end{cases}
\end{equation}
Take $d=T+1+Tm$ and
\[
 c_0=\dfrac{1}{138240},\qquad C_0=12800\sqrt{60},\qquad c_1=\dfrac{27}{320}\cdot\dfrac{0.6}{1280\cdot60^{1/4}\sqrt{24}} \approx2.9\times10^{-6}.
\]

\subsection{Chain certificates}\label{app:nestedimport}
For the nonconvex chain of \citet{cdhs}, write
\[
 \hat f_{\mathrm{nc}}(\bar x)=\dfrac12(1-\bar x_1)^2
 +\dfrac12\sum_{i=1}^T(\bar x_i-\bar x_{i+1})^2
 +\sum_{i=1}^T\Upsilon(\bar x_i).
\]
We use the following facts:
\begin{enumerate}
\item \cite[Lemma~3.4]{cjz}. If $\nu=\theta^2/(3m-1)$, then
\[
 \min_y\hat Q_m(a,y,b)=\nu(a-b)^2,\qquad
 F(\bar x):=\min_{x^{(1)},\ldots,x^{(T)}}\hat f(x)=\nu\hat f_{\mathrm{nc}}(\bar x).
\]
\item \cite[Lemma~3.5]{cjz}. If $\theta\sqrt m\le1$, then $\norm{\nabla\hat f(x)}\ge\norm{\nabla F(\bar x)}/2$.
\item \cite[Lemmas~2 and~3]{cdhs}, restated in \citet[Lemma~3.2]{cjz}.
\[
 \begin{gathered}
 \Upsilon'(0)=\Upsilon'(1)=0,\qquad \Upsilon\ge\Upsilon(1)=0,\qquad \Upsilon(0)\le10,\qquad \sup|\Upsilon''|\le360,\\
 \sup|\Upsilon'''|\le480,\qquad \bar x_T=\bar x_{T+1}=0\quad\Longrightarrow\quad \norm{\nabla\hat f_{\mathrm{nc}}(\bar x)}>\dfrac14.
 \end{gathered}
\]
\end{enumerate}
Direct differentiation gives
\[
 \Upsilon''(t)=120\dfrac{t^4+3t^2-2t}{(1+t^2)^2},\qquad
 \Upsilon'''(t)=240\dfrac{-t^3+3t^2+3t-1}{(1+t^2)^3}.
\]
The bounds follow from
\[
 \begin{gathered}
 t^4\le(1+t^2)^2,\qquad 3t^2\le\dfrac34(1+t^2)^2,\qquad
 2|t|\le1+t^2,\\
 |-t^3+3t^2+3t-1|
 =|\operatorname{Im}(t+\mathrm i)^3-\operatorname{Re}(t+\mathrm i)^3|
 \le\sqrt2(1+t^2)^{3/2},\\
 \Upsilon''(0)=0,\qquad
 \Upsilon(0)=120\int_0^1\dfrac{t^2(1-t)}{1+t^2}\,dt
 \le120\int_0^1t^2(1-t)\,dt=10.
 \end{gathered}
\]
For the gradient transfer, let $H_m=\nabla_y^2(\hat Q_m/2)$ and let $y^*(\bar x)$ minimize the bridge terms. The two endpoint inverse columns have equal norms, with
\[
 (H_m^{-1}e_1)_j=\dfrac{m(2m-j)}{3m-1},\qquad
 \norm{H_m^{-1}e_1}\le\dfrac23m^{3/2},\qquad
 \nabla_y\hat f=H_m(y-y^*).
\]
Each outer coordinate meets at most two endpoints. Hence, for $\theta\sqrt m\le1$,
\[
 \norm{\nabla_{\bar x}\hat f-\nabla F} \le\dfrac43\theta\sqrt m\,\norm{\nabla_y\hat f},\qquad \norm{\nabla F} \le\norm{\nabla_{\bar x}\hat f}+\dfrac43\norm{\nabla_y\hat f} \le\dfrac53\norm{\nabla\hat f}\le2\norm{\nabla\hat f}.
\]

\subsection{Membership and stationarity}
In either accuracy range of Theorem~\ref{thm:nested},
\[
 \begin{gathered}
 \nu\le8\Bigl(\dfrac{c_0}{60}\Bigr)^{1/2}\le\dfrac1{12},\qquad
 360\delta\nu\le360\Bigl(\dfrac{c_0}{60}\Bigr)^{1/2}L=\dfrac L8,\\
 m=\lfloor(3\nu)^{-1/2}\rfloor\ge2,\qquad
 \theta^2m=m(3m-1)\nu\le3m^2\nu\le1.
 \end{gathered}
\]
Since $\inf\hat f_{\mathrm{nc}}=0$, the gap satisfies
\[
 f(0)-\inf f =\delta\nu\beta^2\hat f_{\mathrm{nc}}(0) \le\delta\nu\beta^2(1/2+10T) \qquad \le20T\delta\nu\beta^2 =\dfrac{1280T\eps^2}{\delta\nu}\le\Delta.
\]
In particular,
\begin{equation}\label{eq:Texplicit}
 T=\Bigl\lfloor\dfrac{\rho^{1/2}\Delta\eps^{-3/2}}{1280\sqrt{60}}\Bigr\rfloor\ge10
 \qquad\text{if }\rho^{1/2}\Delta\eps^{-3/2}\ge C_0.
\end{equation}
Only the coordinatewise potentials $\delta\nu\beta^2\Upsilon(\bar x_i/\beta)$ are non-quadratic, so
\[
 \norm{\nabla^3f_k(x)}=\norm{\nabla^3C(x)}
 \le480\delta\nu\beta^{-1}
 =\dfrac{60(\delta\nu)^2}{\eps}=\rho.
\]
The common quadratic Hessian is a sub-sum of the positive semidefinite quadratic Hessian in \citet[Appendix~D]{cjz}. Equal-owner edges are vertex-disjoint, and $\sum_kP_k$ is a union of path Laplacians. Thus
\[
 \norm{\nabla^2C}\le4\delta+360\delta\nu,\qquad
 P_k^2=2P_k,\qquad \norm{P_k}\le2,\qquad
 \Bigl\|\sum_kP_k\Bigr\|\le4.
\]
For IS, $\delta=L/(8n)$ and $n\ge4$ give
\[
 \norm{\nabla^2f_k}\le4\delta+360\delta\nu+2n\delta
 \le\dfrac L8+\dfrac L8+\dfrac L4<L.
\]
For AS, put $u=x-y$ and $g=\nabla C(x)-\nabla C(y)$. Then
\[
 \begin{aligned}
 \dfrac1n\sum_{k=1}^n\norm{\nabla f_k(x)-\nabla f_k(y)}^2 &\le2\norm g^2+2n\delta^2\sum_k\norm{P_ku}^2 =2\norm g^2+4n\delta^2u^\top\Bigl(\sum_kP_k\Bigr)u\\
 &\le\bigl[2(4\delta+360\delta\nu)^2+16n\delta^2\bigr]\norm u^2 \le\Bigl(\dfrac9{32}+\dfrac14\Bigr)L^2\norm u^2<L^2\norm u^2,
 \end{aligned}
\]
using $\delta=L/(8\sqrt n)$, $4\delta\le L/4$, and $360\delta\nu\le L/8$. Membership holds for every owner assignment. Finally,
\[
 \bar x_T=\bar x_{T+1}=0\quad\Longrightarrow\quad
 \norm{\nabla f(x)}=\delta\beta\norm{\nabla\hat f(x/\beta)}
 >\dfrac{\delta\beta\nu}{8}=\eps.
\]

\subsection{Support propagation and owner search}
Index coordinates by \eqref{eq:order}. Every coupling joins adjacent coordinates, and $\Upsilon'(0)=\Upsilon''(0)=0$. Therefore
\[
 \supp(x)\cup\supp(v)\subseteq[p]
 \quad\Longrightarrow\quad
 \supp(\nabla f_k(x))\cup\supp(\nabla^2f_k(x)v)\subseteq[p+1].
\]
If the frontier link is an internal edge $e$, position $p+1$ can be exposed only when $k=\sigma_e$. For the union $R_t$ of reply supports and $p_t=\max R_t$, Definition~\ref{def:zr} gives
\[
 p_0=0,\qquad p_t\le p_{t-1}+1,\qquad
 p_t>p_{t-1}\ \Longrightarrow\ i_t=\sigma_e
 \quad\text{at an internal frontier edge }e.
\]

\begin{lemma}[Owner search]\label{lem:search}
Let $\mathcal E$ be the internal edges of bridges $1,\ldots,T-1$, so $J=\abs{\mathcal E}=(T-1)(m-1)$. For a zero-respecting algorithm, let $Y_e$ be the number of queries made while the frontier link is $e$, with $Y_e=\infty$ if the frontier never reaches $e$. Then $\Pr[\sum_{e\in\mathcal E}Y_e\le(n-1)J/4]\le e^{-J/8}$.
\end{lemma}
\begin{proof}
Enumerate $\mathcal E$ as $e_1,\ldots,e_J$ in \eqref{eq:order}, and let $\mathcal G_r$ contain the transcript and private randomness when the frontier reaches $e_r$. Extend a stopped policy beyond its budget by searching successive frontier owners; its original transcript and output are unchanged. Before reaching $e_r$, all endpoints of $e_r$ and later edges vanish in the query points and directions, so their owners do not affect any reply. Conditionally on $\mathcal G_r$, $\sigma_{e_r}$ is uniform on
\[
 S_r=\begin{cases}
 [n],&e_r\text{ is the first edge of a bridge},\\
 [n]\setminus\{\sigma_{e_{r-1}}\},&\text{otherwise}.
 \end{cases}
\]
At $e_r$, an informative failed query excludes only its queried owner; a point and direction vanishing at the frontier give no information. Hence the owner's rank $U_r$ in any adaptive search order is conditionally uniform on $[|S_r|]$, is $\mathcal G_{r+1}$-measurable, and satisfies
\[
 Y_{e_r}\ge U_r,\qquad |S_r|\in\{n-1,n\},\qquad
 \E[U_r\mid\mathcal G_r]\ge n/2,\qquad 1\le U_r\le n.
\]
Conditional Hoeffding, with $s=Jn/2-(n-1)J/4\ge(n-1)J/4$, gives
\[
 \Pr\Bigl[\sum_rU_r\le\dfrac{(n-1)J}{4}\Bigr]
 \le\exp\Bigl(-\dfrac{2s^2}{J(n-1)^2}\Bigr)\le e^{-J/8}.\qedhere
\]
\end{proof}

\subsection{Completion of the proof}
\begin{proof}
Assume first that $n\ge4$. For a budget $N<(n-1)J/4$, an output with $\bar x_T\ne0$ or $\bar x_{T+1}\ne0$ requires traversing every edge in $\mathcal E$. Since $J\ge9$,
\[
 \Pr[\norm{\nabla f(\hat x)}\le\eps]
 \le\Pr\Bigl[\sum_{e\in\mathcal E}Y_e\le N\Bigr]
 \le e^{-J/8}\le e^{-9/8}<\dfrac13.
\]
Consequently,
\[
 \Nzr_\eps\ge\dfrac{(n-1)(T-1)(m-1)}4
 \ge\dfrac{27}{320}nTm,
\]
using $n-1\ge3n/4$, $T-1\ge9T/10$, and $m-1\ge m/2$. Moreover,
\[
 m\ge\dfrac23(3\nu)^{-1/2},\qquad
 T\ge\dfrac{0.9\rho^{1/2}\Delta\eps^{-3/2}}{1280\sqrt{60}}.
\]
Substitution of \eqref{eq:nestedparams} yields
\[
 nTm\ge\dfrac{0.6}{1280\cdot60^{1/4}\sqrt{24}}
 \begin{cases}
 n^{3/4}L^{1/2}\rho^{1/4}\Delta\eps^{-7/4}&\text{(AS)},\\
 \sqrt n\,L^{1/2}\rho^{1/4}\Delta\eps^{-7/4}&\text{(IS)}.
 \end{cases}
\]

For $n\le3$, take identical components $f_k=\delta\beta^2\hat f(\cdot/\beta)$ with $\delta=L/8$ and the other parameters of \eqref{eq:nestedparams}. The Hessian-Lipschitz, gap, and gradient conditions above depend on $\delta$ only through $\delta\nu$, and $\norm{\nabla^2f_k}\le4\delta+360\delta\nu<L$ since $\nu\le1/90$ in the stated range. A component query returns a full query, and $p_t\le p_{t-1}+1$, so an output with $\bar x_T\ne0$ or $\bar x_{T+1}\ne0$ requires at least $(T-1)(m-1)\ge0.45\,Tm$ queries. With $\nu=8(\rho\eps/60)^{1/2}/L$, this exceeds $3^{3/4}c_1L^{1/2}\rho^{1/4}\Delta\eps^{-7/4}\ge c_1n^{3/4}L^{1/2}\rho^{1/4}\Delta\eps^{-7/4}$, so the same constant $c_1$ serves.
\end{proof}

\section{Proofs for Unrestricted Randomized Algorithms}\label{app:highlower}
\subsection{The replicated instance}\label{app:replicated}
With smoothness budget $L_\star$ and the universal constants of \citet{zhangrand}, set
\begin{equation}\label{eq:zhangparams}
 \begin{gathered}
 h=H_0\eps,\qquad \gamma=c_\gamma(\rho h)^{1/2},\qquad D_0=h^2/\gamma,\\
 M=\lfloor c_M(L_\star/\gamma)^{1/2}\rfloor,\qquad
 N=\lfloor c_N\Delta/D_0\rfloor,\qquad m=\lfloor M/2\rfloor.
 \end{gathered}
\end{equation}
For $\eps\le c_0L_\star^2/\rho$ and $\rho^{1/2}\Delta\eps^{-3/2}\ge C_0$, \citet[Lemmas~5.2--5.4]{zhangrand} give $m,N\ge2$ and
\begin{equation}\label{eq:depth}
 Nm\asymp L_\star^{1/2}\rho^{1/4}\Delta\eps^{-7/4}.
\end{equation}
They also provide distinct eigenvalues $\lambda_j\in[\gamma,L_\star/4]$, positive weights $\nu_j$, and a residual
\[
 r(\lambda)=\lambda^{-1}-p_\star(\lambda),\qquad
 \deg p_\star\le m-1,\qquad \norm r_\nu\ge c/\gamma,
\]
orthogonal to polynomials of degree below $m$ in the $\nu$-weighted spectral inner product. Put $\Lambda=\operatorname{diag}(\lambda_jI_{K_j})$, with multiplicities specified below.

Fix $Nn$ public, mutually orthogonal owner subspaces of dimension $\sum_jK_j$. Independently draw Haar rotations $U_{si}$ and uniform unit vectors $u_{si}$, and define
\begin{equation}\label{eq:replicated}
 A_{si}=U_{si}\Lambda U_{si}^\top,\qquad
 b_{si}=\dfrac h{\sqrt n}U_{si}u_{si},\qquad
 a_s=\dfrac{r(A_s)b_s}{\norm{r(A_s)b_s}},
\end{equation}
where $A_s,b_s$ collect the owner blocks and $z_{si}$ denotes owner $i$'s coordinates in block $s$. Choose an arbitrary unit readout on the null event of a zero denominator. Using $R_B,\vartheta,\psi$ from \citet{zhangrand}, with $\psi^{(k)}=0$ on $(-\infty,1/4]$ for every $k\ge0$, put
\[
 v_s(z_s)=\dfrac{a_s^\top R_B(z_s)}{\alpha_0},\qquad
 \ell_s(z_s)=D_0[C_\ell+\vartheta(b_s^\top z_s/D_0)],
\]
\begin{equation}\label{eq:repinstance}
 G(z) =-\ell_1(z_1)-\sum_{s=2}^N\psi\bigl(v_{s-1}(z_{s-1})\bigr)\ell_s(z_s),\qquad f_i(z) =G(z)+\dfrac n2\sum_{s=1}^Nz_{si}^\top A_{si}z_{si}.
\end{equation}
The average
\[
 F(z)=G(z)+\dfrac12\sum_sz_s^\top A_sz_s
\]
uses the same cascade as \citet{zhangrand}. Each $A_s$ has the distinct eigenvalues $\lambda_j$ with multiplicities $nK_j$. Writing $\Pi_j$ for the latent projector onto the $\lambda_j$ eigenspace of one owner, its normalized spectral weights are
\[
 \omega_{sj}=\dfrac1n\sum_{i=1}^n\norm{\Pi_ju_{si}}^2.
\]
Thus the aggregate weights lie within the same relative error of $\nu_j$ whenever every owner's weights do.

\subsection{Regularity and spectral certificates}\label{app:zhangimport}
The norm-based bounds of \citet[Proposition~6.1]{zhangrand} hold for every realization:
\[
 F(0)-\inf F\le\Delta,\qquad
 \norm{\nabla^2F}\le L_\star,\qquad
 \Lip(\nabla^2F)\le\rho.
\]
The owner Hessians $Q_i$ have orthogonal supports, $\norm{Q_i}\le L_\star/4$, and $\norm{\nabla^2G}\le5L_\star/4$. Hence
\[
 \begin{gathered}
 \norm{\nabla^2f_i}\le\dfrac{n+5}{4}L_\star,\qquad \dfrac1n\sum_i\norm{\nabla f_i(x)-\nabla f_i(x')}^2 \le\dfrac{n+25}{8}L_\star^2\norm{x-x'}^2,\\
 L_\star=\begin{cases}L/(2\sqrt n)&\text{(AS)},\\L/(2n)&\text{(IS)}.\end{cases}
 \end{gathered}
\]
Since $\nabla^3f_i=\nabla^3F$, every component is individually $\rho$-Hessian-Lipschitz.

The stationarity and projection certificates of \citet[Lemmas~6.3 and~7.2]{zhangrand} require
\[
 \norm{b_s}=h,\qquad \gamma I\preceq A_s\preceq\dfrac{L_\star}{4}I,
\]
and spectral weights within $1\pm\delta$ of $\nu$. The first two hold exactly. If $J$ is the number of distinct eigenvalues, choose, for the horizon $T$ of Lemma~\ref{lem:span},
\[
 K_j \ge C\bigl[(2T)^2+\eta^{-2}(2T)+\delta^{-2}\bigr] \log\dfrac{2JNn}{\zeta},\qquad \left|\dfrac{K_j}{\sum_kK_k}-\nu_j\right| \le\dfrac{\delta\nu_j}{4}.
\]
The multiplicities are chosen after fixing a horizon $T\ge4nNM+1$, including one slot for the terminal output. All owner-weight and Gaussian projection bounds then hold with probability at least $1-\zeta$. The norm and moment arguments of \citet[Lemma~6.3]{zhangrand} depend only on these spectral weights and the displayed norm bounds, so they apply to the replicated block. On this event,
\[
 \norm{\nabla F(z)}\le\eps\quad\Longrightarrow\quad
 \begin{cases}
 \psi(v_s(z_s))=1&1\le s<N,\\
 \norm{A_Nz_N-b_N}\le\eps,\\
 \norm{z_s}\le B&1\le s\le N.
 \end{cases}
\]
In particular $R_B(z_N)=z_N$. The moment bound $a_N^\top A_N^{-1}b_N\ge\alpha_0$ and the fixed choices of $c_\alpha,H_0$ give
\[
 a_N^\top z_N\ge a_N^\top A_N^{-1}b_N-\eps/\gamma
 \ge3\alpha_0/4,\qquad v_N(z_N)\ge3/4.
\]

\subsection{Scheduled disclosure}
Set $\alpha_0=c_\alpha h/\gamma$ and $B=C_Bh/\gamma$ as in \citet{zhangrand}. Choose universal $c_{\mathrm{delay}}>0$ small enough, and the accuracy constant so that $m\ge\lceil2/c_{\mathrm{delay}}\rceil$. Then
\[
 q_0=\lfloor c_{\mathrm{delay}}nm\rfloor
 \ge\dfrac12c_{\mathrm{delay}}nm\ge1.
\]
Before query selection, the auxiliary oracle discloses $(u_{si},b_{si})_i$ at $\tau_s=1+(s-1)q_0$ and $(U_{si})_i$ at $\sigma_s=1+sq_0$. Each call also reveals the latent coordinates of its point and direction in every owner subspace, and their products with $A_{si}$ for the selected owner. Its censored reply is the value, gradient, and Hessian--vector product of
\[
 \begin{gathered}
 k(t) =\min\{N,1+\lfloor(t-1)/q_0\rfloor\},\qquad G_t(z) =-\ell_1(z_1)-\sum_{s=2}^{k(t)}\psi\bigl(v_{s-1}(z_{s-1})\bigr)\ell_s(z_s),\\
 f_{i,t}(z) =G_t(z)+\dfrac n2\sum_{s=1}^Nz_{si}^\top A_{si}z_{si}.
 \end{gathered}
\]
Every retained coupling has a disclosed predecessor rotation and successor forcing vector. We show below that omitted couplings are flat at all query points.

\subsection{Selected-owner information span}\label{app:span}
\begin{lemma}\label{lem:span}
Fix a deterministic policy and a horizon $T\ge4nNM+1$, and let the owner dimension be large enough. For every pair $(s,i)$ there are $2T$ jointly independent standard Gaussian vectors, independent of $u_{si}$, such that for $t\le T$ and $t<\sigma_s$ both latent query vectors of the $t$th call, restricted to $(s,i)$, lie in
\[
 V_{si}+\mathbf 1\{t\ge\tau_s\}\spn\{u_{si},\Lambda u_{si},\ldots,\Lambda^{q_{si}(t)}u_{si}\},
\]
where $V_{si}$ is the $\Lambda$-cyclic span of the Gaussian vectors and $q_{si}(t)$ counts calls to component $i$ in $[\tau_s,t)$. The Gaussian and independence assertions hold also conditional on the raw draws of every other pair. The same conclusion holds for a passive point selected from the preceding transcript and exposed without returning a matrix product or incrementing an owner counter.
\end{lemma}
\begin{proof}
Fix $(s,i)$, condition on the raw draws of all other pairs, and write $U=U_{si}$, $u=u_{si}$. Before $\sigma_s$, the retained couplings involve this pair only through its disclosed forcing vector: its outgoing readout is censored, and every retained predecessor readout has a disclosed rotation. The censored replies are therefore computable from the other draws, the disclosed forcing vector, the latent point and direction, and the selected-owner matrix products, which reduces the policy to an adaptive strategy for the single pair. Conditioning on raw independent draws rather than on adaptive transcripts preserves their distribution.

Let $\mathcal L$ be the span of the latent vectors whose images under $U$ have been exposed: earlier latent points and directions, their products with $\Lambda$ on calls to owner $i$, and $u$ once disclosed. Conditional on this history, the restriction $U:\mathcal L^\perp\to(U\mathcal L)^\perp$ is Haar. Before $\tau_s$ the history does not use $u$, so this also holds conditional on $u$; at the fixed round $\tau_s$, disclosing $(u,b_{si})$ adjoins $u$ to $\mathcal L$ and preserves the conditional Haar law on the complement.

At round $t$, expose the point first, without returning a matrix product. Decompose its physical restriction $w=w_\parallel+w_\perp$ relative to $U\mathcal L$. If $w_\perp\ne0$, then $e=U^\top w_\perp/\norm{w_\perp}$ is conditionally uniform on the unit sphere of $\mathcal L^\perp$, also conditional on $u$. Draw independently $a\sim N(0,P_{\mathcal L})$ and a chi radius $R$ with $\dim\mathcal L^\perp$ degrees of freedom, and put $g=a+Re$. Then $g$ is conditionally $N(0,I)$ and $U^\top w\in\mathcal L+\spn\{g\}$; if $w_\perp=0$, a fresh standard Gaussian serves the same purpose. Adjoin $U^\top w$ to $\mathcal L$ and treat the direction, which was fixed before the point exposure, in the same way. Each exposure reveals $U$ only on the new direction $e$, whose image is the fixed vector $w_\perp/\norm{w_\perp}$, and the added Gaussian component and radius are independent. Hence the conditional Haar invariant persists, and the successive conditional laws make all innovations jointly independent and independent of $u$.

Just before call $t$, let $q=q_{si}(t)$ and assume inductively that all exposed latent vectors lie in $V_{si}+\mathbf 1\{t\ge\tau_s\}\spn\{u,\Lambda u,\ldots,\Lambda^qu\}$. The two exposures add one Gaussian start each and preserve the inclusion. Only afterwards are the products returned. If component $i$ was selected, applying $\Lambda$ to both latent vectors raises the marked degree to at most $q+1$, not twice, since neither query vector depends on the current products; the counter increases by one. Otherwise no product for this pair is returned and the counter is unchanged. Gaussian terms stay in $V_{si}$, which proves the inclusion. A passive exposure is the same argument with no product and no increment. At most $4T+1$ generators enter $\mathcal L$; pad the Gaussian family with independent draws to length $2T$. A sufficiently large owner dimension, ensured by the multiplicity choice, completes the construction.
\end{proof}

\subsection{Proof of Theorem~\ref{thm:highlower}}
\begin{proof}
On the spectral-weight event, write
\[
 a_{si}=\beta_{si}U_{si}\tilde a_{si},\qquad
 \tilde a_{si}=\dfrac{r(\Lambda)u_{si}}{\norm{r(\Lambda)u_{si}}},\qquad
 \sum_i\beta_{si}^2=1,\qquad \beta_{si}^2\le\dfrac2n
\]
for $\delta\le1/3$. At $t<\sigma_s$, define the heavy owners by $\mathcal H_s(t)=\{i:q_{si}(t)\ge m\}$. Then
\[
 |\mathcal H_s(t)|\le\dfrac{q_0}{m}\le c_{\mathrm{delay}}n,\qquad
 \sum_{i\in\mathcal H_s(t)}\beta_{si}^2\le2c_{\mathrm{delay}}.
\]
The mass bound holds for every subset, including this adaptive choice. For each light owner's latent query vector $x$, Lemma~\ref{lem:span} and \citet[Lemma~7.2]{zhangrand} give
\[
 |\tilde a_{si}^\top x|\le C(\delta+\eta)\norm x.
\]
Radiality of $R_B$ and Cauchy--Schwarz across owners yield
\begin{equation}\label{eq:readout}
 |a_s^\top R_B(z_s)|
 \le2B\left(\sqrt{2c_{\mathrm{delay}}}+C(\delta+\eta)\right)
 \le\dfrac{\alpha_0}{4},
\end{equation}
for sufficiently small universal $c_{\mathrm{delay}},\delta,\eta$, since $B/\alpha_0=C_B/c_\alpha$ is universal.

Every omitted coupling has a flat gate, so its value and all derivatives vanish. In particular its Hessian--vector contribution is zero for every queried direction. Induction on rounds couples the censored and true transcripts on the concentration event. For $t<\tau_N=\sigma_{N-1}$, the $(N-1)$st gate is closed, contradicting the stationarity certificate; for $\tau_N\le t<\sigma_N$, \eqref{eq:readout} contradicts its consequence $v_N(z_N)\ge3/4$.

To include an arbitrary terminal output, fix a query cap $Q$, pad earlier stopping with ignored zero queries, and freeze the output after reply $Q$, before any disclosure scheduled for round $Q+1$. At that round, expose its latent coordinates as a passive point with zero direction, no matrix product, and no counter increment. The output is measurable with respect to the enlarged history, so Lemma~\ref{lem:span} and \eqref{eq:readout} apply to it. If $c\le c_{\mathrm{delay}}/4$, then $Q<c\,nNm$ implies $Q+1\le Nq_0<\sigma_N$, by $q_0\ge c_{\mathrm{delay}}nm/2$ and $N\ge2$, so the output also violates a necessary stationarity condition.

With $\zeta=1/6$, every deterministic policy with fewer than $c\,nNm$ calls fails with probability at least $5/6$ under the hard distribution, whose multiplicities depend only on the horizon and the class parameters. Fixing the private randomness of a randomized algorithm and integrating gives the same failure probability, so no algorithm succeeds with probability $2/3$ on every instance.

Finally, \eqref{eq:depth} gives
\[
 nNm\asymp
 \begin{cases}
 n^{3/4}L^{1/2}\rho^{1/4}\Delta\eps^{-7/4}&L_\star=L/(2\sqrt n),\\
 \sqrt n\,L^{1/2}\rho^{1/4}\Delta\eps^{-7/4}&L_\star=L/(2n),
 \end{cases}
\]
with the accuracy thresholds
\[
 \dfrac{c_0L_\star^2}{\rho}
 =\begin{cases}
 c_0L^2/(4n\rho)&\text{(AS)},\\
 c_0L^2/(4n^2\rho)&\text{(IS)}.
 \end{cases}
 \qedhere
\]
\end{proof}

\section{Proofs for the Low-Accuracy Lower Bound}\label{app:lowlower}
\subsection{The dense-weak-hiding instance}\label{app:dwh}
For $j\in[T]$, take $x_j\in\R^D$, a hidden direction $\Theta^{(j)}\in\{\pm1\}^D$, and a fixed encoding table $S^{(j)}\in\{\pm1\}^{n\times D}$ satisfying
\[
 \dfrac1n\sum_iS^{(j)}_{i:}=\varrho\,\Theta^{(j)},\qquad
 n^{-1/2}\le\varrho<\dfrac{17}{16}n^{-1/2}.
\]
For $\Psi(z)=z/(1+\norm z^2)^{1/2}$, define
\[
 q_j(z)=\Bigl\langle\dfrac{\Theta^{(j)}}{\sqrt D},\Psi(z)\Bigr\rangle,\qquad
 \tilde q_{ij}(z)=\dfrac1\varrho\Bigl\langle\dfrac{S^{(j)}_{i:}}{\sqrt D},\Psi(z)\Bigr\rangle,
\]
\begin{equation}\label{eq:ptj}
 R_i(x)=-\tilde q_{i1}(x_1)-\sum_{j=2}^{T}G\bigl(q_{j-1}(x_{j-1})\bigr)\bigl(\tilde q_{ij}(x_j)+1\bigr)+\dfrac{33}{32}\norm{x}^2,
\end{equation}
where $G$ is the gate of \citet[Lemma~A.1]{densehiding}, with all derivatives zero on $(-\infty,5/32]$. The scaling $f_i(y)=\alpha R_i(\beta y)$ of \citet[Equation~(148)]{densehiding} gives individual smoothness $L$ and gradient norm at least $8\eps$ before the last stage is revealed. Recovering each hidden direction $\Theta^{(j)}$ requires order $\min\{n,\varrho^{-2}\}=n$ distinct rows.

\paragraph{The average is free of the bias.}
Since $n^{-1}\sum_i\tilde q_{ij}=q_j$, the average
\begin{equation}\label{eq:ptjavg}
 R(x)=-q_1(x_1)-\sum_{j=2}^{T}G\bigl(q_{j-1}(x_{j-1})\bigr)\bigl(q_j(x_j)+1\bigr)+\dfrac{33}{32}\norm{x}^2
\end{equation}
is independent of $\varrho$, whereas component derivatives scale as $1/\varrho$. Appendix~\ref{app:third} and the scaling relations
\[
 \sup_x\norm{\nabla^3R(x)}\le\ell_3\le1.23\times10^9,\qquad
 \alpha\beta^2=\dfrac{L\varrho}{341},\qquad
 \alpha\beta=\dfrac{8\eps}{\chi},
\]
where $\chi=0.279247\ldots$ is from \citet[Lemma~A.6]{densehiding}, give
\begin{equation}\label{eq:avghl}
 \Lip\Bigl(\nabla^2\dfrac1n\sum_if_i\Bigr) \le\alpha\beta^3\ell_3 =\dfrac{(\alpha\beta^2)^2}{\alpha\beta}\ell_3 \le c_H\dfrac{L^2}{n\eps},\qquad c_H =\dfrac{(17/16)^2\chi\ell_3}{8\cdot341^2}<414.
\end{equation}

\subsection{The third derivative of the average}\label{app:third}
\begin{lemma}\label{lem:Q3}
For a unit vector $u$, the function $q_u(x)=\langle u,x\rangle(1+\norm x^2)^{-1/2}$ satisfies $\sup_x\norm{D^3q_u(x)}=3$, with equality at $x=0$.
\end{lemma}
\begin{proof}
For a unit vector $h$, put $S=1+\norm x^2$ and $c=\langle x,h\rangle$. Differentiation gives $D^3q_u(x)[h,h,h]=\langle u,V\rangle$, where
\[
 V=x\bigl(9cS^{-5/2}-15c^3S^{-7/2}\bigr)+h\bigl(-3S^{-3/2}+9c^2S^{-5/2}\bigr).
\]
Writing $x=ch+w$, $w\perp h$, $z=c^2/S$, and $e=\norm w^2/S$, we have $z,e\ge0$, $z+e<1$, and
\[
 \dfrac{\norm V^2}{9} =(1-z-e)^3\bigl[(-1+6z-5z^2)^2+ze(3-5z)^2\bigr] \le(1-z)^4(5z^2-2z+1)\le1.
\]
The last function is nonincreasing on $[0,1]$, with derivative $-2(1-z)^3(15z^2-10z+3)$. Banach's theorem for symmetric multilinear forms identifies the norm with the supremum over diagonal unit inputs.
\end{proof}

\begin{lemma}[Third derivative of the dense-weak-hiding average]\label{lem:third}
The average \eqref{eq:ptjavg} satisfies $\sup_x\norm{\nabla^3R(x)}\le\ell_3$ with $\ell_3\le1.2212\times10^9$, independently of $n$, $\varrho$, $T$, and $D$.
\end{lemma}
\begin{proof}
Put $Q_k=\sup_{x,u}\norm{\nabla^kq_u(x)}$ and $M_k=\norm{G^{(k)}}_\infty$. By \cite[Lemma~A.4]{densehiding}, Lemma~\ref{lem:Q3}, and \cite[Lemma~A.1]{densehiding},
\[
 Q_1\le1,\quad Q_2\le2,\quad Q_3=3,\qquad M_1\le17,\quad M_2\le259.
\]
For $G(s)=H_0(8s-5/4)$, its derivative is a convolution of a triangular density with the normalized bump $\phi_g$ of half-width $d_g=1/512$ \cite[Equations~(43)--(46)]{densehiding}. The density's distributional second derivative has point masses with weights $4/w_g^2$, $-8/w_g^2$, $4/w_g^2$, where $w_g=255/256$. The translated bumps are disjoint, hence
\[
 M_3=512\cdot8w_g^{-2}\norm{\phi_g}_\infty.
\]
For $a_g=(1-16d_g^2)^{-1}$, normalization yields
\[
 \norm{\phi_g}_\infty^{-1} =d_g^2\int_{-1/d_g}^{1/d_g} \exp\!\left(-\dfrac{t^2}{1-d_g^2t^2}\right)dt \ge d_g^2\int_{-4}^{4}e^{-a_gt^2}\,dt \ge d_g^2\left(\sqrt{\dfrac\pi{a_g}} -\dfrac{e^{-16a_g}}{4a_g}\right),
\]
using the Gaussian-tail bound $2\int_4^\infty e^{-a_gt^2}dt\le e^{-16a_g}/(4a_g)$. Thus $\norm{\phi_g}_\infty<147904$ and $M_3<6.10576\times10^8$.

For $A_j=g(x_{j-1})p(x_j)$, with $g=G\circ q_{j-1}$ and $p=q_j+1\in(0,2)$, put $r_j=\norm{h_j}$ and
\[
 \gamma_1=M_1Q_1,\qquad
 \gamma_2=M_2Q_1^2+M_1Q_2,\qquad
 \gamma_3=M_3Q_1^3+3M_2Q_1Q_2+M_1Q_3.
\]
Then
\[
 |D^3A_j[h,h,h]| \le2\gamma_3r_{j-1}^3+3\gamma_2Q_1r_{j-1}^2r_j +3\gamma_1Q_2r_{j-1}r_j^2+Q_3r_j^3.
\]
Using $a^2b\le(2a^3+b^3)/3$ and $ab^2\le(a^3+2b^3)/3$, and counting each $r_k$ in at most two links, gives
\[
 |D^3R[h,h,h]|\le\ell_3\sum_kr_k^3\le\ell_3,\qquad
 \ell_3=2\gamma_3+3\gamma_2Q_1+3\gamma_1Q_2+Q_3.
\]
This constant includes the contribution $Q_3r_1^3$ of $-q_1$; the quadratic contributes zero. Banach's theorem converts the diagonal bound to the operator norm.
\end{proof}

\begin{lemma}[Finite-order replies]\label{lem:jets}
The conclusion of \citet[Theorem~3.1]{densehiding} holds when every call returns all derivatives of $f_i$ at the query point up to any fixed finite order.
\end{lemma}
\begin{proof}
Conditions (R2)--(R3) of \citet[Corollary~B.6]{densehiding} require non-completing replies to depend only on the history and queried rows of the current stage, with no contribution from unopened links. Since $G^{(r)}(s)=0$ for $r\ge0$ and $s\le5/32$, all derivatives of unopened links vanish at query points whose omitted gate inputs lie below the threshold. Completion precedes the reply, so Lemma~B.4, Theorem~B.5, and Corollary~B.6 of that work apply. This concerns finite-order derivatives at the query point, not local germs: a flat gate need not vanish on a neighborhood of its threshold. The auxiliary $\Omega(n)$ instance \cite[Appendix~C.5]{densehiding} is quadratic.
\end{proof}

\begin{proposition}[General form of the low-accuracy lower bound]\label{prop:lowgeneral}
There are universal constants $n_0,c_d>0$ such that, for $n\ge n_0$ and $\eps^2\le c_d\Delta L_{\mathrm{eff}}$ with $L_{\mathrm{eff}}=\min\{L,(n\rho\eps/c_H)^{1/2}\}$ and $c_H$ as in \eqref{eq:avghl}, every randomized algorithm needs
\[
 \Omega\!\left(n+\min\left\{\sqrt n\,L\Delta\eps^{-2},\,n\rho^{1/2}\Delta\eps^{-3/2}\right\}\right)
\]
queries on $\Fis(n,L,\rho,\Delta)$, even when calls return full component Hessians. Below the crossover $\eps<L^2/(n\rho)$, the initialization condition of Theorem~\ref{thm:matching} ensures the gap condition for a sufficiently large universal $C_0$. For the terms proportional to $\Delta$, the ratio of the FastCubic upper rate to the displayed full-Hessian lower rate is, up to logarithmic factors, of order $(n\rho\eps)^{1/2}/L$ when $\eps\ge L^2/(n\rho)$ and of order $(L^2/(n\rho\eps))^{1/4}$ when $\eps\le L^2/(n\rho)$. The additive $n$ term can reduce the ratio of the total query bounds.
\end{proposition}
\subsection{Proof of the lower bound in Theorem~\ref{thm:low}}\label{app:lowproof}
\begin{proof}
Apply \citet[Theorem~3.1]{densehiding} and Lemma~\ref{lem:jets} with $L_{\mathrm{eff}}=\min\{L,(n\rho\eps/c_H)^{1/2}\}$. The instance is individually $L_{\mathrm{eff}}$-smooth, and \eqref{eq:avghl} gives $\Lip(\nabla^2f)\le c_HL_{\mathrm{eff}}^2/(n\eps)\le\rho$. For $n\ge n_0$ and $\eps^2\le c_d\Delta L_{\mathrm{eff}}$, the resulting lower bound is
\[
 \Omega\!\left(n+\sqrt n\,L_{\mathrm{eff}}\Delta\eps^{-2}\right)
 =\Omega\!\left(n+\min\left\{
 \sqrt n\,L\Delta\eps^{-2},\,
 n\rho^{1/2}\Delta\eps^{-3/2}\right\}\right).
\]
In the low-accuracy regime, $L_{\mathrm{eff}}\ge L/\sqrt{c_H}$, so $c_3\le c_d/\sqrt{c_H}$ ensures the gap condition. Below the crossover, it is equivalent to
\[
 \rho^{1/2}\Delta\eps^{-3/2}\ge\dfrac{\sqrt{c_H}}{c_d\sqrt n},
\]
which follows from Theorem~\ref{thm:matching}'s initialization condition for sufficiently large universal $C_0$. Both source instances also satisfy averaged smoothness \cite[Corollary~3.2]{densehiding}. This proves Proposition~\ref{prop:lowgeneral} and the lower bound in Theorem~\ref{thm:low}.
\end{proof}

\section{Proofs for Section~\ref{sec:results}}\label{app:matching}
\subsection{Proofs of Theorems~\ref{thm:matching} and~\ref{thm:is}}
\begin{proof}
Theorems~\ref{thm:nested} and~\ref{thm:highlower} give the lower bounds. For the upper bound, apply Theorem~\ref{thm:fastcubic} and use
\[
 n\rho^{1/2}\Delta\eps^{-3/2}
 \le c_0^{1/4}n^{3/4}L^{1/2}\rho^{1/4}\Delta\eps^{-7/4},
 \qquad \eps\le c_0L^2/(n\rho).
\]
The condition $\rho^{1/2}\Delta\eps^{-3/2}\ge C_0$ absorbs initialization and one cubic subproblem. The individual-smoothness upper bound follows from $\Fis\subseteq\Fas$.
\end{proof}

\subsection{Proof of Proposition~\ref{prop:hessian}}
\begin{proof}
Proposition~\ref{prop:lowgeneral} gives the lower bound. Cubic regularization~\cite{nesterov} and PathVR give $O(n+n\rho^{1/2}\Delta\eps^{-3/2})$ and $\tO(n+\sqrt n\,L\Delta\eps^{-2})$, respectively. Taking the smaller bound proves the claim.
\end{proof}

\needspace{8\baselineskip}
\subsection{Proof of Proposition~\ref{prop:single}}\label{app:single}
\begin{proof}
Take $n=1$, $\tau_s=1+(s-1)m$, and $\sigma_s=1+sm$. By Lemma~\ref{lem:span}, before $\sigma_s$ the latent queries lie in $V_s+\spn\{u_s,\ldots,\Lambda^{m-1}u_s\}$. The projection estimate keeps every predecessor gate flat, so no point queried before $1+Nm$ is $\eps$-stationary. As in the proof of Theorem~\ref{thm:highlower}, the terminal output is one passive exposure at round $Q+1<1+Nm$ for a budget $Q<Nm$, so it obeys the same projection bound. Equation~\eqref{eq:depth} with $L_\star=L$ gives $\Omega(Nm)=\Omega(L^{1/2}\rho^{1/4}\Delta\eps^{-7/4})$.
\end{proof}

\end{document}